\documentclass{paper}
\usepackage{url}
\usepackage{amsthm}
\usepackage{amssymb,amsmath,mdwlist, mathtools}
\usepackage[linesnumbered,ruled,vlined]{algorithm2e}
\usepackage{subcaption}
\usepackage{authblk}
\usepackage{xcolor}
\title{A sequential approach for speed planning under jerk constraints}

\author{Luca~Consolini,
        Marco~Locatelli,
        and~Andrea~Minari
\thanks{The authors are all with Dipartimento di Ingegneria e
  Architettura, Universit\`a di Parma\\Parco Area delle Scienze,
  181/A, Parma, Italy, e-mail: luca.consolini@unipr.it,
  marco.locatelli@unipr.it, andrea.minari2@studenti.unipr.it}}

\newtheorem{defn}{Definition}
\newtheorem{thm}{Theorem}
\newtheorem{remark}{Remark}
\newtheorem{obser}{Observation}
\newtheorem{prop}{Proposition}
\newtheorem{lem}{Lemma}

\newtheorem{problem}{Problem}

\newcommand{\Real}{\mathbb{R}}

\newcommand{\bq}{\mathbf{q}}

\newcommand{\bd}{\mathbf{d}}

\newcommand{\bgamma}{\boldsymbol{\gamma}}

\newcommand{\etab}{\boldsymbol{\eta}}

\newcommand{\dw}{\delta w}
\newcommand{\dwb}{\boldsymbol{\delta}\mathbf{w}}

\newcommand{\bA}{b_{\text{A}}}
\newcommand{\bdd}{b_{\text{D}}}
\newcommand{\bn}{b_{\text{N}}}
\newcommand{\bp}{b_{\text{P}}}

\newcommand{\bAi}{b_{\text{A}_i}}
\newcommand{\bddi}{b_{\text{D} i}}
\newcommand{\bni}{b_{\text{N}_i}}
\newcommand{\bpi}{b_{\text{P}_i}}

\newcommand{\ub}{{\mathbf u}}
\newcommand{\xb}{{\mathbf x}}
\newcommand{\yb}{{\mathbf y}}

\newcommand{\by}{{\mathbf y}}
\newcommand{\wb}{{\mathbf w}}
\newcommand{\nub}{\boldsymbol{\nu}}

\newcommand{\scalar}[2]{\left<#1,#2\right>}

\providecommand{\keywords}[1]
{
  \small	
  \textbf{\textit{Keywords---}} #1
}

\begin{document}

\maketitle

\begin{abstract}
In this paper we discuss a sequential algorithm for the computation of a minimum-time speed profile over a given path, under velocity, acceleration and jerk constraints. 
Such a problem arises in industrial contexts such as automated warehouses, where LGVs need to perform assigned tasks 
as fast as possible in order to increase productivity.
It can be reformulated as an optimization problem with a convex objective function, linear velocity and acceleration constraints, and non-convex jerk constraints, which, thus, represent the main source of difficulty. While existing non-linear programming (NLP) solvers can be employed for the solution of this problem, it turns out that the performance and robustness of such solvers can be enhanced by the sequential line-search algorithm proposed in this paper. At each iteration a feasible direction, with respect to the current feasible solution, is computed, and a step along such direction is taken in order to compute the next iterate. The computation of the feasible direction is based on the solution of a linearized version of the problem, and the solution of the linearized problem, through an approach which strongly exploits its special structure, represents the main contribution of this work. The efficiency of the proposed approach with respect to existing NLP solvers is proved through different computational experiments.  
\end{abstract}


\section{Introduction}

An important problem in motion planning is the computation of the
minimum-time motion of a car-like vehicle from a start configuration
to a target one while avoiding collisions (obstacle avoidance) and
satisfying kinematic, dynamic, and mechanical constraints (for
instance, on velocities, accelerations and maximal steering
angle).
This problem can be approached in two ways: i) As a
minimum-time trajectory planning, where both the path to be
followed by the vehicle and the timing law on this path (i.e., the
vehicle's velocity) are simultaneously designed. For instance, one
could use the RRT algorithm (see~\cite{7505628}).
ii) As a (geometric) path planning followed by a minimum-time
  speed planning on the planned path (see, for
  instance,~\cite{doi:10.1177/027836498600500304}).
 In this paper, following the second paradigm, we assume that the path
that joins the initial and the final configuration is assigned, and we
aim at finding the time-optimal speed law
that satisfies some kinematic and dynamic constraints.
The problem can be reformulated as an optimization one and
it is quite relevant from the practical point of view. In particular, 
in automated warehouses the speed of LGVs needs to be planned under acceleration and jerk constraints. The solution algorithm should be: i) {\em fast}, since speed planning is made continuously throughout the work-day, not only when an LGV receives a new task but also during the execution of the task itself, since conditions may change, e.g., if the LGV has to be halted for security reasons; ii)  {\em reliable}, i.e., it should return solutions of high quality, because a better speed profile
allows to save time and even a small percentage improvements, say a 5\% improvement, has a considerable impact on the productivity of the warehouse and, thus, determines a significant economic gain.
In our previous work~\cite{consolini2017scl}, we proposed an optimal
time-complexity
algorithm for finding the time-optimal speed law that satisfies
constraints on maximum velocity and tangential and normal acceleration.
In the subsequent work~\cite{CabConLoc2018coap1}, we included a bound
on the derivative of the acceleration with respect to the arc-length.
In this paper, we consider the
presence of jerk constraints (constraints on the time
derivative of the acceleration). The resulting optimization problem
is a non-convex one and, for this reason, is significantly more
complex than the ones we discussed in~\cite{consolini2017scl}
and~\cite{CabConLoc2018coap1}. The main contribution of this work is
the development of a line-search algorithm for this problem based on
the sequential solution of convex problems.
The proposed algorithm meets both the requirement of being fast and the requirement of being reliable. The former is met by heavily exploiting the special structure of the optimization problem, the latter by the theoretical guarantee that
the returned solution is a first-order stationary point (in practice, a minimizer) of the optimization problem.


\subsection*{Problem Statement}
\label{sec:velplan}
Here we introduce more formally the problem at hand.
Let $\boldsymbol{\gamma}: [0,s_f] \to \mathbb{R}^2$ be a smooth
function. The image set $\boldsymbol{\gamma}([0,s_f])$ is the path to be
followed, $\boldsymbol{\gamma}(0)$ the initial
configuration, and $\boldsymbol{\gamma}(s_f)$ the final one. 
Function $\boldsymbol{\gamma}$ has arc-length parameterization, that is, is such that
\mbox{$(\forall \lambda \in [0,s_f])$,  $\|\boldsymbol{\gamma}'(\lambda)\| =1$}. 
In this way, $s_f$ is the length of the path.
We want to compute the speed-law that
minimizes the overall transfer time (i.e., the time needed to go from
$\boldsymbol{\gamma}(0)$ to $\boldsymbol{\gamma}(s_f)$).
To this end, let $\lambda: [0,t_f] \to
[0,s_f]$ be a differentiable monotone increasing function, that
represents the vehicle's curvilinear abscissa as a function of time, and
let
$v: [0,s_f]\to [0,+\infty[$ be such that $(\forall t \in
[0,t_f])\ \dot \lambda(t)=v(\lambda(t))$. In this way, $v(s)$
is the derivative of the vehicle curvilinear abscissa, which corresponds to the norm of its velocity vector at position $s$. The position of the vehicle as a function of time is given by ${\bf x}:[0,t_f] \to \Real^2, \ {\bf x}(t)=\boldsymbol{\gamma}(\lambda(t))$.
The velocity
and acceleration are given, respectively, by
$$
\begin{array}{ll}
\dot {\bf x}(t)=\boldsymbol{\gamma}'(\lambda(t)) v(\lambda(t)),\\
\ddot {\bf x}(t)=a_T(t) \boldsymbol{\gamma}'(\lambda(t))+ a_N(t) \boldsymbol{\gamma}'^{\perp} (\lambda(t)),\,
\end{array}
  $$
where $a_T(t)= v'(\lambda(t)) v(\lambda(t))$, $a_N(t)=k(\lambda(t))
v(\lambda(t))^2$ are, respectively, the tangential and normal
components of the acceleration (i.e., the projections of the acceleration
vector $\ddot {\bf x}$ on the tangent and the normal to
the curve). Moreover $\boldsymbol{\gamma}'^{\perp}(\lambda)$ is the normal to
vector $\boldsymbol{\gamma}'(\lambda)$, the tangent of
$\boldsymbol{\gamma}'$ at $\lambda$.
Here $k:[0,s_f] \to \Real$ is the
scalar curvature, defined as
\mbox{$k(s)=\scalar{\boldsymbol{\gamma}''(s)}{\boldsymbol{\gamma}'(s)^\perp}$}.
Note that $|k(s)|=\|\boldsymbol{\gamma}''(s)\|$.
In the following, we
assume that $k(s) \in \mathcal{C}^{1}([0,s_f],\Real)$.
The total maneuver time, for a given velocity profile \mbox{$v\in C^1([0,s_f],\Real)$}, is returned by the functional
\begin{equation}
\label{obj_fun_pr}
{\cal F}:  C^1([0,s_f],\Real)\rightarrow \Real, \ \ \ {\cal F}(v)=\int_0^{s_f} v^{-1}(s) d s.
\end{equation}
In our previous work~\cite{consolini2017scl}, we considered the
problem
\begin{equation}
\label{eqn_problem_pr}
\min_{v \in {\cal V}}  {\cal F}(v),
\end{equation} 
where the feasible region ${\cal V}\subset  C^1([0,s_f],\Real)$ is defined by the following set of constraints
\begin{subequations}
\label{eqn_problem_constraints}
\begin{align}
v(0)=0,\,v(s_f)=0, \label{inter_con_pr}\\
0\leq v(s) \leq  v_{\max},\ \  s \in ]0,s_f[, \label{con_speed_pr}\\
|2 v'(s)v(s)| \leq A, \ \ s \in [0,s_f],  \label{con_at_pr}\\
|k(s)| v(s)^2 \leq A_N, \ \ s \in [0,s_f],  \label{con_an_pr}
\end{align}
\end{subequations}
($v_{\max}$, $A$, $A_N$ are upper bounds for the velocity, the tangential acceleration, and the normal acceleration, respectively).
Constraints~\eqref{inter_con_pr} are the initial and final
interpolation
conditions, while constraints~\eqref{con_speed_pr},~\eqref{con_at_pr},~\eqref{con_an_pr}
limit velocity and the tangential and normal components of
acceleration.  In~\cite{consolini2017scl} we presented an algorithm, with linear-time computational
complexity with respect to the number of variables, that provides an
optimal solution of~\eqref{eqn_problem_pr} after spatial
discretization. One limitation of this algorithm is that the
obtained velocity profile is Lipschitz\footnote{A function $f:\Real
  \to \Real$ is \emph{Lipschitz} if there exists a real positive
  constant $L$ such that $(\forall x,y \in \Real)\ |f(x)-f(y)| \leq L |x-y|$.
  } but not
differentiable, so that the vehicle's acceleration is
discontinuous. With the aim of obtaining a smoother
velocity profile, in the subsequent work~\cite{CabConLoc2018coap1}, we required that the velocity be differentiable and
we imposed a Lipschitz condition (with constant $J$) on its derivative.
In this way, after setting $w=v^2$, the feasible region of the problem ${\cal W}\subset
  C^1([0,s_f],\Real)$ is defined by the set of functions $w \in
  C^1([0,s_f],\Real)$ that satisfy the following set of constraints
\begin{subequations}
\label{eqn_problem_cont_constraints}
\begin{align}
w(0)=0,\,w(s_f)=0, \label{inter_con_cont}\\
0\leq w(s) \leq  v_{\max}^2,\ \ s \in ]0,s_f[, \label{con_speed_cont}\\
|w'(s)| \leq A,\ \ s \in [0,s_f],  \label{con_at_cont}\\
|k(s)| w(s) \leq A_N, \ \ s \in [0,s_f],  \label{con_an_cont}\\
|w'(s_1)- w'(s_2)| \leq J |s_1-s_2|,\ \  s_1,s_2 \in [0,s_f]. \label{con_j_cont}
\end{align}
\end{subequations}
So that we end up with
problem
\begin{equation}
\label{eqn_problem_cont}
\min_{w\in {\cal W}}  G(w),
\end{equation} 
where the objective function is
\begin{equation}
\label{obj_fun_cont}
G:  C^1([0,s_f],\Real)\rightarrow \Real, \ \ \ G(w)=\int_0^{s_f} w^{-1/2}(s) d s.
\end{equation}
The objective function~\eqref{obj_fun_cont} and constraints
\eqref{inter_con_cont}-\eqref{con_an_cont}
correspond to the ones in Problem~\eqref{eqn_problem_pr} after the
substitution $w=v^2$. Note that this change of variable is well known in the literature. It has been first proposed in \cite{Pfeiffer87}, while
in \cite{verscheure09} it is observed that Problem \eqref{eqn_problem_pr} becomes convex after this change of variable.
The added set of constraints~\eqref{con_j_cont} is a Lipschitz
condition on the derivative of the squared velocity $w$. It is used to
enforce a smoother velocity profile by bounding the second derivative of the
squared velocity with respect to arc-length. Note that
constraints~\eqref{eqn_problem_cont_constraints} are linear
and that objective function~\eqref{obj_fun_cont} is convex.
In~\cite{CabConLoc2018coap1}, we proposed an algorithm for solving a
finite dimensional approximation of Problem~\eqref{eqn_problem_cont_constraints}.
The algorithm exploited the particular structure of the resulting
convex finite dimensional problem.
This paper extends the results of~\cite{CabConLoc2018coap1}. It
considers a non-convex variation of Problem~\eqref{eqn_problem_cont_constraints}, in which
constraint~\eqref{con_j_cont} is
substituted with a constraint on the time derivative of the
acceleration $|\dot a(t)| \leq J$, where $a(t)=\frac{d}{d t}
v(\lambda(t))=v' (\lambda(t)) v (\lambda(t))=\frac{1}{2} w'(\lambda(t))$.

\begin{figure}[!h]
	\centering
	\includegraphics[width=0.99\linewidth]{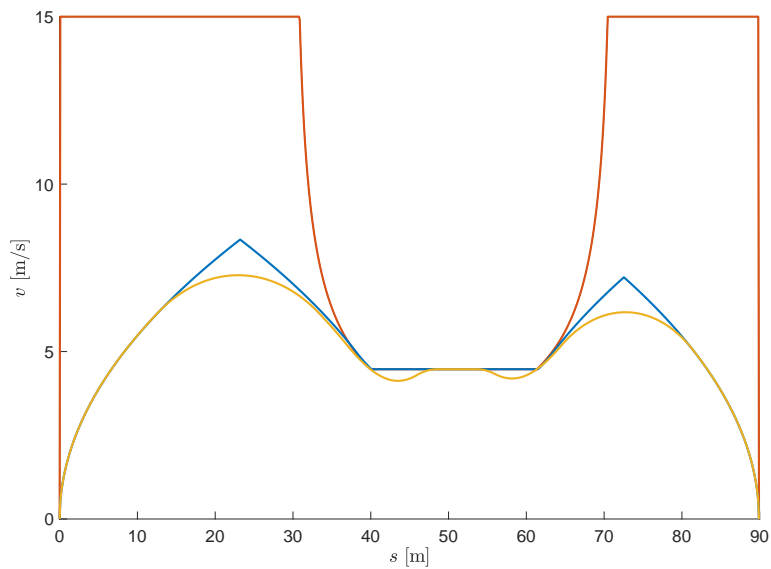}
	\caption{
		Red line represents the maximum allowed velocity along the path as a function of the arc length. The blue line is the optimal speed profile without jerk constraints. The orange line is the optimal speed profile with jerk constraints.
		\label{fig:minaripathjerk}}
\end{figure}

\begin{figure}
	\centering
	\includegraphics[width=0.99\columnwidth]{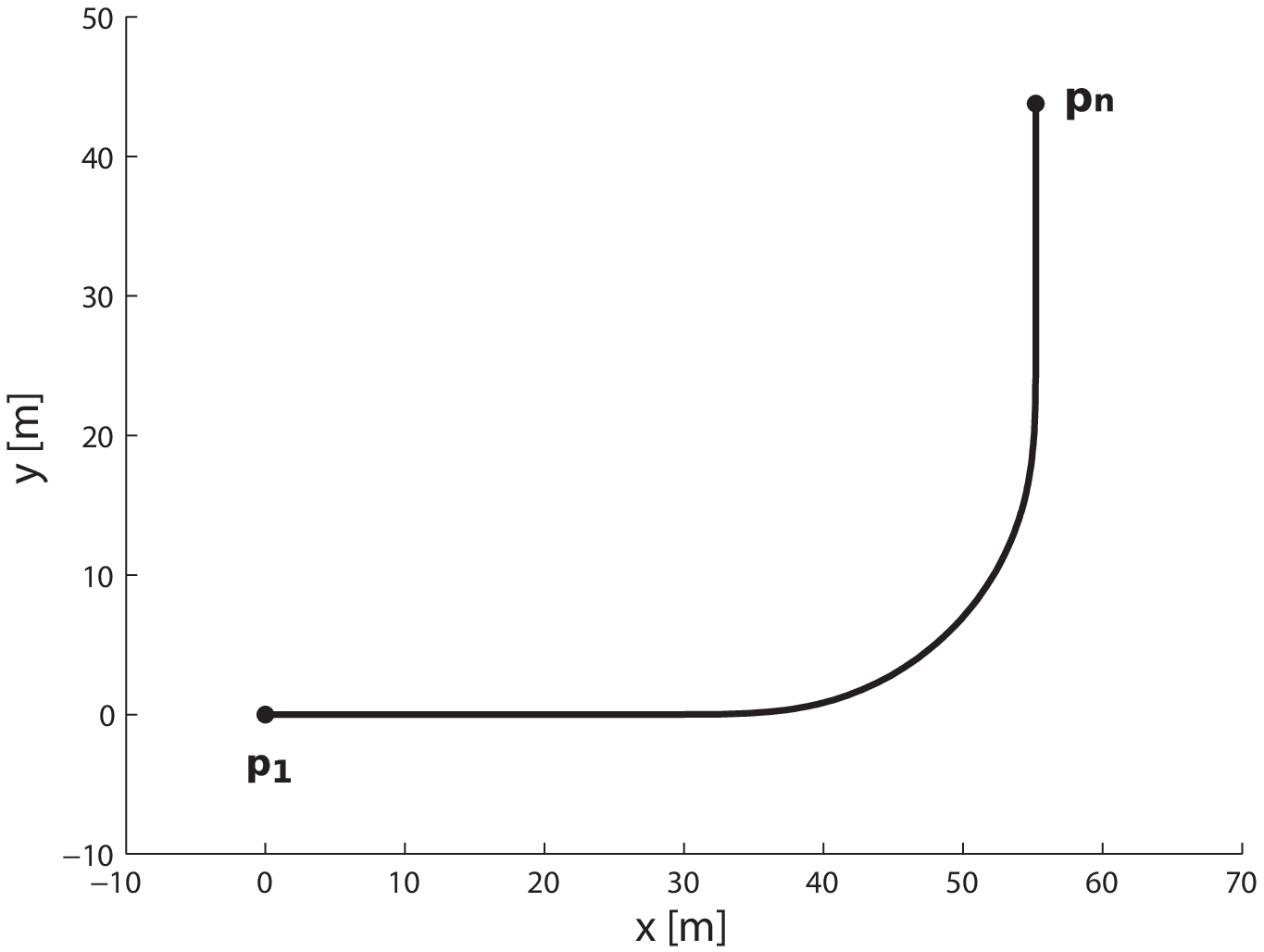}
	\caption{Path with parameters  $s_{f} = 90$ m, maximum allowed velocity  $v_{\max}=15 ms^{-1}$, $A = 1.5$ ms$^{-2}$, $A_N = 1$ ms$^{-2}$, $J=1$ ms$^{-3}$}
	\label{cap1:fig:geometry_path}
\end{figure}

Then, we set
$$j_L(t) = \dot{a}(t) = \frac{1}{2}
w^{\prime\prime}(s(t))\sqrt{(w(s(t)))}.$$
We name this quantity ``jerk''. Note  that $j_L(t)$ is not the
third time derivative of the position $\bgamma(\lambda(t))$ but is
related
to it. We will clarify this in Section~\ref{sec_higher_dim}.
We used the subscript $L$ for ``longitudinal''.
Then, we end up with the following minimum-time problem:
\begin{problem}[Smooth minimum-time velocity planning problem: continuous version]
	\label{cap4:prob:continuos}
	\begin{align}
	&  \min_{w \in C^2} {\displaystyle\int_0^{s_f} w(s)^{-1/2} \, ds} \nonumber\\
	&w(0)=0, \quad w(s_{f})=0 , \nonumber \\\
	&0  \leq  w(s) \leq \mu^+(s), & s \in [0,s_{f}],  \nonumber\\
	&\frac{1}{2}\left|w^{\prime}(s)\right| \leq A, & s\in [0,s_{f}], \label{cap4:acc_bound_cont}\\
	&\frac{1}{2}\left|w^{\prime\prime}(s) \sqrt{w(s)}\right| \le J,  & s\in[0,s_f],\label{cap4:jerk_bound_cont}
	\end{align}
\end{problem}
where $\mu^{+}$ is the square velocity upper bound depending on the shape of the path, i.e.,
\[
\mu^{+}(s) = \min	\left\{v_{\max}^2,\frac{A_N}{|k(s)|}\right\},
\] 
with $v_{\max}$, $A_N$ and $k$ be the maximum allowed velocity of the vehicle, the maximum normal acceleration and the curvature of the path, respectively.
Parameters $A$ and $J$
are the bounds representing the limitations on the (tangential) acceleration and the jerk, respectively.
For the sake of simplicity we consider constraints~\eqref{cap4:acc_bound_cont} 
and~\eqref{cap4:jerk_bound_cont} symmetric and constant. However, the following development could be easily extended 
to the non-symmetric and non-constant case. Note that the jerk constraint~\eqref{cap4:jerk_bound_cont} is non-convex.
The continuous problem is discretized as follows. We subdivide the path into $n-1$ intervals of equal length, i.e.,
we evaluate function $w$ at points
$$
s_i=\frac{(i-1) s_f}{n-1},\ \ \ i=1,\ldots,n,
$$
so that we have the following $n$-dimensional vector of variables
$$
{\bf w}=(w_1, w_2,\ldots,w_n)=\left(w(s_1), w(s_2), \ldots, w(s_n)\right).
$$
Then, the finite dimensional version of the problem is: 
\begin{problem} [Smooth minimum-time velocity planning problem: discretized version]
	\label{cap4:prob_disc}
	\begin{align}
	& \qquad \min_{\wb \in \Real^n} \sum_{i=1}^{n-1} \frac{2h}{\sqrt{w_{i+1}} + \sqrt{w_{i}}}
	\label{cap4:obj:disc}\\
	&0  \leq  \wb \leq \ub, \label{cap4:con:bound}\\ 
	& w_{i+1} - w_i \leq 2hA, & i=1, \dots ,n-1,\label{cap4:con:acc}\\ 
	&w_{i} - w_{i+1} \leq 2hA, & i=1, \dots ,n-1, \,\label{cap4:con:dec} \\
	&(w_{i-1} - 2w_i + w_{i+1})\sqrt{\frac{w_{i+1} + w_{i-1}}{2}}\le 2h^2J,& i = 2,\dots,n-1, \label{cap4:con:par}\\
	&-(w_{i-1} - 2w_i + w_{i+1})\sqrt{\frac{w_{i+1} + w_{i-1}}{2}}\le 2h^2J,& i = 2,\dots,n-1, \label{cap4:con:nar}
	\end{align}
\end{problem}
where $u_i=\mu^+(s_i)$, for $i=1,\ldots,n$, and, in particular, $u_1 = 0$ and $u_n = 0$, since we are assuming that the initial and final velocity are equal to 0.  
The objective function (\ref{cap4:obj:disc}) is an approximation
of~\eqref{obj_fun_cont} given by the Riemann sum of the intervals obtained by dividing each interval $[s_i,s_{i+1}]$, for
$i=1,\ldots,n-1$, in two subintervals of the same size. 
Constraints (\ref{cap4:con:acc}) and (\ref{cap4:con:dec}) are obtained by finite difference to approximate $w^{\prime}$.
Constraints~\eqref{cap4:con:par} and~\eqref{cap4:con:nar} 
are obtained by using a second-order central finite
difference to approximate $w^{\prime\prime}$, while $w$ is approximated by  the arithmetic mean.
Due to jerk constraints~\eqref{cap4:con:par} and~\eqref{cap4:con:nar},
Problem~\ref{cap4:prob_disc} is a non-convex one and
cannot be solved with the algorithm presented
in~\cite{CabConLoc2018coap1}. 
\newline\newline\noindent
For the sake of illustration, Figure~\ref{fig:minaripathjerk} shows the speed profiles computed with and without the jerk constraints
for the instance associated to the path shown in Figure~\ref{cap1:fig:geometry_path}. It is straightforward to observe how the jerk bounded velocity profile is smoother than the one obtained without the jerk limitation.
It is also interesting to remark that in the interval where the maximum allowed velocity is smallest (the one between 40 and 60 m), the optimal solution falls below the maximum speed profile at the beginning and at  the end of the interval. Indeed, due to the jerk constraints,
it is worthwhile to reduce the speed in these positions in order to have larger velocities in the two intervals immediately preceding  and immediately following it, respectively. 
\subsection*{Main result}
The main contribution of this paper is the development of a new
solution algorithm for finding a local minimum of the non-convex
Problem~\ref{cap4:prob_disc}.
As will be detailed in next sections, we propose to solve Problem~\ref{cap4:prob_disc} by a  line-search algorithm based on the sequential solution of convex problems. 
The algorithm will be an iterative one where at each iteration the following operations will be performed.
\newline\newline\noindent
{\bf Constraint linearization:} We will first define a convex problem by linearizing constraints~\eqref{cap4:con:par} and~\eqref{cap4:con:nar} through a first-order Taylor approximation around the current point 
$\wb^{(k)}$. Differently from other sequential algorithms for non-linear programming (NLP) problems, we will keep the original convex objective function. The linearized problem will be introduced in Section \ref{sec:linesearch}.
\newline\newline\noindent
{\bf Computation of a feasible descent direction:} The convex problem (actually, a relaxation of such problem) is solved in order to compute a feasible descent direction $\dwb^{(k)}$. The main contribution of the paper lies in this part. 
The computation requires the minimization of a suitably defined objective function through a further iterative algorithm. At each iteration of this algorithm the following operations are performed:
\begin{itemize}
\item {\bf Objective function evaluation:} Such evaluation requires the solution of a problem with the same objective function but subject to a subset of the constraints. The special structure of the resulting subproblem is heavily exploited in order to solve it efficiently.
This is the topic of Section  \ref{cap4:sec:accnar}. \newline\newline\noindent
\item {\bf Computation of a descent step:} Some Lagrange multipliers of the subproblem define a subgradient for the objective function. This can be employed to define a linear programming (LP) problem which returns a descent step for the objective function. This is the topic
of Section~\ref{cap4:sec:PAR}. 
\end{itemize} 
$\ $\newline\newline\noindent
{\bf Line-search:} Finally, a standard line-search along the half-line $\wb^{(k)}+\alpha \dwb^{(k)}$, $\alpha\geq 0$, is performed.
\newline\newline\noindent
Sections \ref{sec:linesearch}--\ref{cap4:sec:PAR} will detail all what we discussed above. In Section \ref{sec_higher_dim}
we will briefly discuss the speed planning problem for a curve in
a generic configuration space and show that, also in this case, a speed profile obtained by
solving Problem~\ref{cap4:prob:continuos} allows to bound the
velocity, the acceleration and the jerk of the obtained trajectory. Finally, in Section \ref{sec:compexp} we will present different computational experiments.
\subsection*{Comparison with existing literature}
Although many works consider the problem of minimum-time speed
planning
with acceleration constraints (see for instance, \cite{7515141,Velenis2008,frego2017}), relatively few consider jerk
constraints.
Perhaps, this is also due to the fact that the jerk constraint is
non-convex, so that its presence significantly increases the complexity of the
optimization task. One can use a general purpose NLP solver (such as SNOPT
or IPOPT) for finding a
local solution of Problem~\ref{cap4:prob_disc}, but the required
time is in general too large for the speed planning application. 
As outlined in the previous subsection, in this work we tackle this problem through an approach based on the solution of a sequence of convex subproblems. There are different approaches in the literature based on the sequential solution of convex subproblems. In~\cite{hauser2014fast} it is first observed that the problem with acceleration constraints but no jerk constraints for robotic manipulators can be reformulated as a convex one with linear constraints, and it is solved by a sequence of LP problems obtained by linearizing the objective function at the current point, i.e., the objective function is replaced by its supporting hyperplane at the current point, and by introducing a trust region centered at the current point.
In \cite{pham2018new,CLMNV19} it is further observed that this problem can be solved very efficiently through the solution of a sequence of 2D LP problems. In \cite{LippBoyd2014} an interior point barrier method is used to solve the same problem based on Newton’s method. Each Newton step requires the solution of a KKT system and an efficient way to solve such systems is proposed in that work. Moving to approaches also dealing with jerk constraints, we mention \cite{debrouwere2013time}. 
In this work it is observed that jerk constraints are non-convex but can be written as the difference of two convex functions. Based on this observation, the authors solve the problem by a sequence of convex subproblems obtained by linearizing at the current point the concave part of the jerk constraints and by adding a proximal term in the objective function which plays the same role as a trust region, preventing from taking too large steps. In \cite{Singh15} a slightly different objective function is considered. Rather than minimizing the travelling time along the given path, the integral of the squared difference between the maximum velocity profile and the computed velocity profile is minimized. After representing
time varying control inputs as products of parametric exponential
and a polynomial functions, the authors reformulate the problem in such a way that its objective function is convex quadratic, while non-convexity lies in difference-of-convex functions. The resulting problem is tackled through the solution of a sequence of convex subproblems obtained by linearizing the concave part of the non-convex constraints. 
In \cite{Palleschi19} the problem of speed planning for robotic manipulators with jerk constraints is reformulated in such a way that non-convexity lies in simple bilinear terms. Such bilinear terms are replaced by the corresponding convex and concave envelopes, obtaining the so called McCormick relaxation, which is the tightest possible convex relaxation of the non-convex problem.
Other approaches dealing with jerk constraints do not rely on the solution of convex subproblems. For instance, in \cite{MacFarlane03} 
a concatenation of
fifth-order polynomials is employed to provide smooth trajectories, which results in quadratic jerk profiles, while in \cite{Haschke08} cubic polynomials are employed, resulting in piecewise constant jerk profiles. The decision process involves the choice of the phase durations, i.e., of the intervals over which a given polynomial applies.
A very recent and interesting approach to the problem with jerk
constraints is~\cite{pham2017structure}. In this work an approach
based on numerical integration is discussed.
Numerical integration has been first applied under acceleration constraints
in \cite{Bobrow85,Shin85}. In~\cite{pham2017structure} jerk
constraints are taken into account. The algorithm detects a position
$s$ along the trajectory where the jerk constraint is singular, that
is, the jerk term disappears from one of the
constraints. Then, it computes
the speed profile up to $s$ by computing two maximum jerk profiles and then
connecting them by a minimum jerk profile, found by a shooting method. In
general, the overall solution is composed of a sequence of various
maximum and minimum jerk profiles.
This approach does not guarantee reaching a local
minimum of the traversal time. Moreover, since
Problem~\ref{eqn_problem_cont_constraints} has velocity and
acceleration constraints, the jerk constraint is singular for all
values of $s$, so that the algorithm presented
in~\cite{pham2017structure} cannot be directly applied to Problem~\ref{eqn_problem_cont_constraints}.

Some algorithms use heuristics to quickly find suboptimal solutions
of acceptable quality.
For instance, \cite{Villagra-et-al2012} proposes an algorithm that applies to curves composed of clothoids, circles and
straight lines. The algorithm does not guarantee local
optimality of the solution. Reference~\cite{RaiCGL2019jerk} presents a
very efficient heuristic algorithm. Also this method does not guarantee global nor local optimality.
Various works in literature consider jerk bounds in the speed optimization problem for
robotic manipulators instead of mobile vehicles. This is a slightly
different problem, but mathematically equivalent to Problem~\eqref{cap4:prob:continuos}.
In particular, paper~\cite{DONG20071941} presents a method based on the solution of a large number of non-linear and non-convex
subproblems. The resulting algorithm is slow, due to the large number
of subproblems; moreover, the authors do not prove its convergence.
Reference~\cite{ZHANG2012472} proposes a similar method that gives
a continuous-time solution. Again, the method is computationally slow, since
it is based on the numerical solution of a large number of
differential equations; moreover, the paper does not contain a proof of
convergence or of local optimality. 
Some other works replace the jerk constraint
with \emph{pseudo-jerk}, that is the derivative of the acceleration
with respect to arc-length, obtaining a constraint analogous
to~\eqref{con_j_cont} and ending up with a convex optimization problem.
For instance,~\cite{8569414} adds to the objective function a pseudo-jerk penalizing
term. This approach is computationally convenient but
substituting~\eqref{cap4:jerk_bound_cont} with~\eqref{con_j_cont} may be overly
restrictive at low speeds.

\subsection*{Statement of contribution}
The method presented in this paper is a sequential convex one which aims at finding a local optimizer of
Problem~\ref{cap4:prob_disc}. To be more precise, as usual with
non-convex problems, only convergence to a stationary point can
usually be proved. However, the fact that the sequence of generated
feasible points is decreasing with respect to the objective function
values usually guarantees that the stationary point is a local minimizer, except in rather pathological cases (see, e.g., \cite[Page 19]{Fletcher:00}).
To our knowledge and as detailed in the
following, this algorithm is more efficient than the ones existing in
literature, since it leverages the special structure of the subproblems obtained as local
approximations of Problem~\ref{cap4:prob_disc}. We discussed this
class of problems in our previous work~\cite{ConLocLAu19Graph}.
This structure
allows computing very efficiently a feasible descent direction
for the main line-search algorithm; it is one of the key elements that allows
us to outperform generic NLP solvers.
\section{A sequential algorithm based on constraint~linearization}
\label{sec:linesearch}
\begin{figure}[!h]
	\centering
	\includegraphics[width=\linewidth]{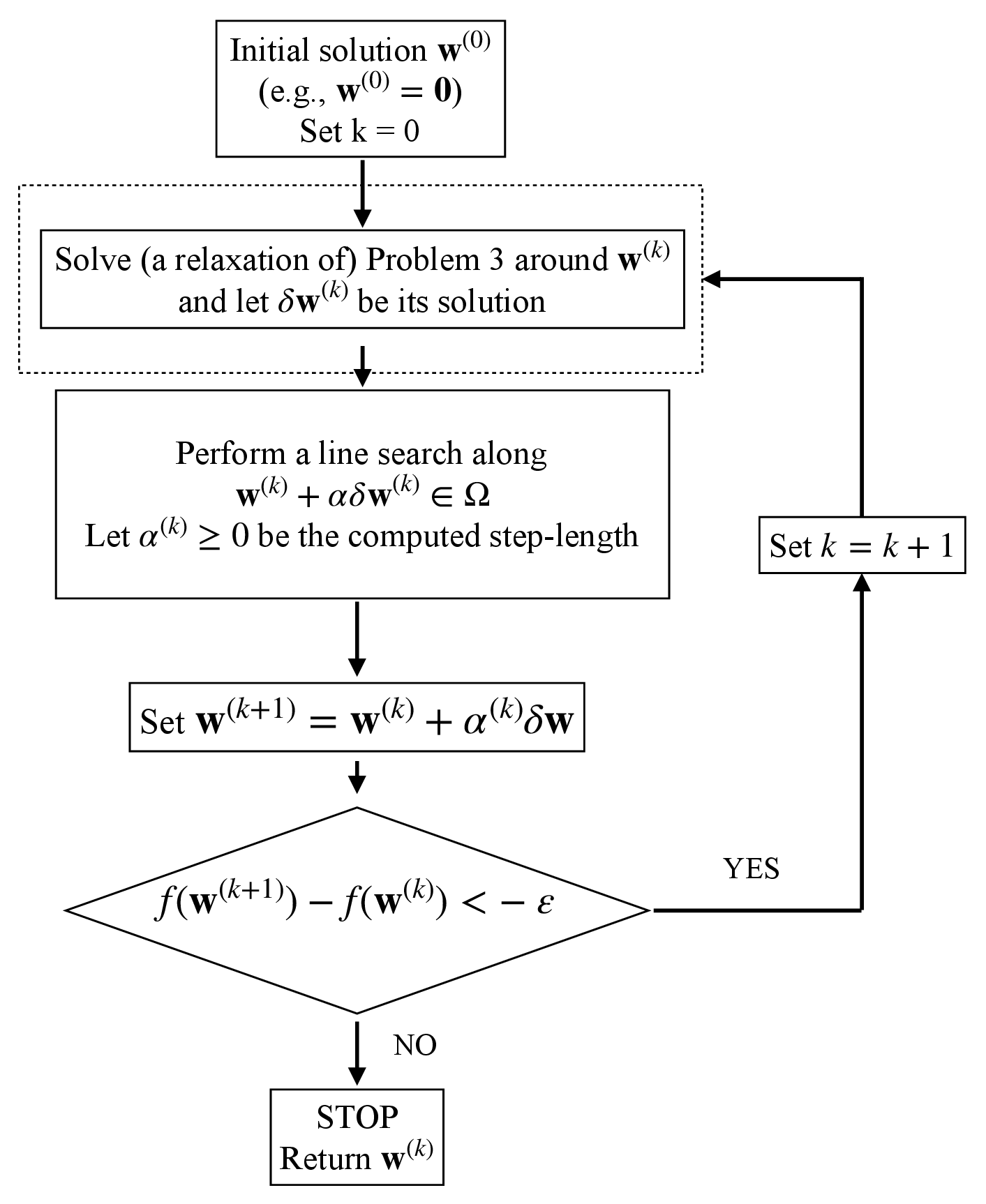}
	\caption{Flow chart of  Algorithm SCA.
The dashed block corresponds to a call of the procedure \texttt{ComputeUpdate}, proposed to solve Problem~\ref{cap4:prob_lin}, which represents the main contribution of this paper.	\label{fig:flowlinesearch}}
\end{figure}

To account for the non-convexity of Problem~\ref{cap4:prob_disc} 
we propose a line-search method 
based on the solution of a sequence of special structured convex problems. Throughout the paper we will call this Algorithm SCA (Sequential Convex Algorithm) and
its flow chart is shown in~Figure~\ref{fig:flowlinesearch}. It belongs to the class of Sequential Convex Programming algorithms, where at each iteration a convex subproblem is solved.
In what follows we will denote by $\Omega$ the feasible region of Problem~\ref{cap4:prob_disc}. 
At each  iteration $k$, we replace the current point $\wb^{(k)}\in\Omega$  with a new point $\wb^{(k)} + \alpha^{(k)}\dwb^{(k)} \in\Omega$, where the step-size $\alpha^{(k)}\in[0,1] $ is obtained by a \emph{line search} along the descent direction $\dwb^{(k)}$, which, in turn, is obtained through the solution of  a convex problem. The constraints of the convex problem  are
linear approximations of  \eqref{cap4:con:bound}-\eqref{cap4:con:nar} around $\wb^{(k)}$,
while the  objective function is the original one. 
Then, the problem we consider to compute the direction $\dwb^{(k)}$ is the following (superscript $k$ of $\wb^{(k)}$ is omitted):
\begin{problem}
	\label{cap4:prob_lin}
	\begin{align}
	& \qquad \min_{\dwb \in \Real^n} \sum_{i=1}^{n-1} \frac{2h}{\sqrt{w_{i+1} +\dw_{i+1}} + \sqrt{w_{i}+\dw_{i}}}
	\label{cap4:eq:obj_lin}\\
	& \mathbf{l_B}  \leq  \dwb  \leq \mathbf{u_B},\label{cap4:eq:bound_lin} \\
	& \dw_{i+1}-\dw_{i}  \le \bAi, & i=1,\dots,n-1,\label{cap4:eq:acc_lin}\\
	&\dw_{i}-\dw_{i+1}  \le \bddi, & i=1,\dots n-1,\label{cap4:eq:dec_lin}\\
	&\dw_{i} - \eta_i \dw_{i-1}  - \eta_i \dw_{i+1} \le \bni ,& i = 2,\dots,n-1.\label{cap4:eq:nar_lin}\\
	& \eta_i \dw_{i-1}  + \eta_i \dw_{i+1} -\dw_{i} \le \bpi ,& i = 2,\dots,n-1,\label{cap4:eq:par_lin}
	\end{align}
\end{problem} 
where $\mathbf{l_B} = -\wb$ and $\mathbf{u_B} = \ub - \wb$ (recall that $\ub$ has been introduced in (\ref{cap4:con:bound}) and its components have been defined immediately below Problem~\ref{cap4:prob_disc}),
while parameters $\etab$, $\mathbf{\bA}$ $\mathbf{\bdd}$, $\mathbf{\bn}$ and $\mathbf{\bp}$ depend on the point $\wb$ around which 
the constraints~\eqref{cap4:con:bound}-\eqref{cap4:con:nar} are linearized.
More precisely, we have:
\begin{equation}
\label{eq:paramdef}
\begin{array}{l}
\bAi = 2hA - w_{i+1} + w_{i} \\ [4pt]
{\bddi} = 2hA - w_{i} + w_{i+1} \\ [4pt]
\eta_i = \frac{3(w_{i+1} + w_{i-1}) -2w_i}{4(w_{i+1}+w_{i-1})} \\ [4pt]
\bpi =  \frac{2\sqrt{2}h^2J - (w_{i-1} - 2w_i + w_{i+1})\sqrt{w_{i+1} + w_{i-1}}}{2\sqrt{w_{i+1}+w_{i-1}}} \\ [4pt]
\bni=  \frac{2\sqrt{2}h^2J + (w_{i-1} - 2w_i + w_{i+1})\sqrt{w_{i+1} + w_{i-1}}}{2\sqrt{w_{i+1}+w_{i-1}}}.
\end{array}
\end{equation}
These parameters are proved to be nonnegative in the following proposition.
\begin{prop}\label{cap4:prop:positivity}
	All  parameters $\etab$, $\mathbf{\bA}$ $\mathbf{\bdd}$, $\mathbf{\bn}$ and $\mathbf{\bp}$ are non negative for $h \rightarrow 0$.
\end{prop}
\begin{proof}
We have $\bAi = 2hA - w_{i+1} + w_{i} \ge 0 $ because of the feasibility of $\wb$.
Analogously we can prove that ${\bddi} = 2hA - w_{i} + w_{i+1} \ge 0 $.
Next, by continuity of $w$ we have that, for $h\to0$, $\eta_i \to \frac{1}{2}$, while by feasibility of $\wb$ we have $\bpi ,\bni \ge 0$.
\end{proof}
The proposed approach follows some standard ideas  of sequential quadratic approaches employed in the literature about non-linearly constrained problems. But a quite relevant difference is that the true objective function~\eqref{cap4:obj:disc} is employed in the problem to compute the direction, rather than a quadratic approximation of such function. This choice comes from the fact that the objective function~\eqref{cap4:obj:disc} has some features (in particular,
convexity and being decreasing), which, combined with the structure of the linearized constraints, allow for an efficient solution of Problem~\ref{cap4:prob_lin}.
Problem~\ref{cap4:prob_lin} is a convex problem with a non-empty feasible region ($\dwb = \mathbf{0}$ is always a feasible solution) and, consequently, can
be solved by existing NLP solvers. 
However, such solvers tend to increase computing times since they need to be called many times
within the iterative Algorithm SCA. 
The main contribution of this paper lies in the routine~\texttt{computeUpdate} (see dashed block in Figure~\ref{fig:flowlinesearch}),
which is able to solve Problem~\ref{cap4:prob_lin} and efficiently returns a descent direction $\dwb^{(k)}$.
To be more precise, we will solve a \emph{relaxation} of Problem~\ref{cap4:prob_lin}. Such relaxation as well as 
the routine to solve it, will be detailed in Sections~\ref{cap4:sec:accnar} and \ref{cap4:sec:PAR}. 
In Section~\ref{cap4:sec:accnar} we present efficient approaches to solve some subproblems including proper subsets of the constraints. Then, in Section~\ref{cap4:sec:PAR} we address the solution of the relaxation of Problem~\ref{cap4:prob_lin}.
\begin{remark}\label{cap4:rem:step-feasible}
	It is possible to see that if one of the
	constraints~\eqref{cap4:con:par}-\eqref{cap4:con:nar} is active at $\wb^{(k)}$, then along the direction $\dwb^{(k)}$ computed through the solution of the linearized Problem \ref{cap4:prob_lin}, it holds
	that $\wb^{(k)} + \alpha \dwb^{(k)} \in\Omega \ $  for any sufficiently small $\alpha>0$. In other words, small perturbations of the current solution $\wb^{(k)}$ along direction $\dwb^{(k)}$ do not lead outside the feasible region
	$\Omega$.
	This fact is illustrated in Figure~\ref{fig:linearization}.
	Let us rewrite constraints~\eqref{cap4:con:par}-\eqref{cap4:con:nar} as follows:
	\begin{equation}\label{cap4:eq:2dconstr}
	|(x-2y)\sqrt{x}|\le C, 
	\end{equation}
	where $x = w_{i+1} + w_{i-1}$, $y = w_i$ and $C=2\sqrt{2}h^2 J$ is a constant. The feasible region associated to constraint~\eqref{cap4:eq:2dconstr} is reported in Figure~\ref{fig:linearization}. In particular, it is the region between the blue and the red curves. 
	\begin{figure}[!h]
		\centering
		\includegraphics[width=0.99\linewidth]{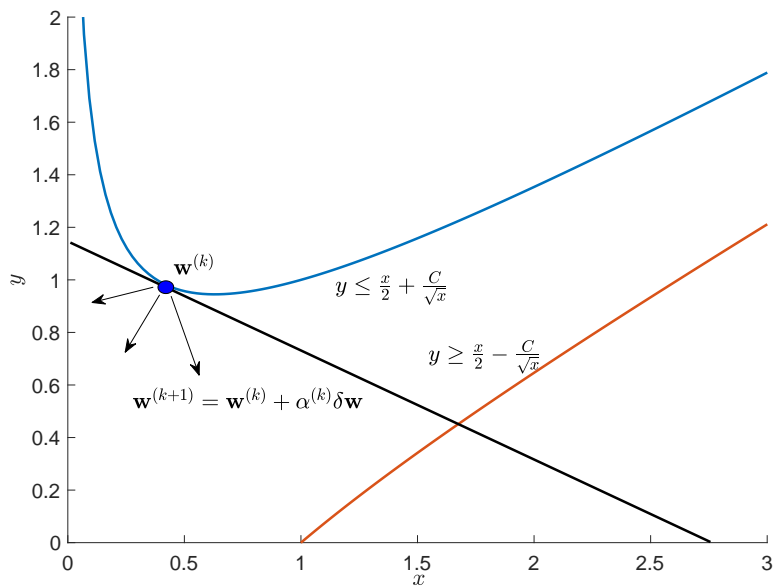}
		\caption{Constraints~\eqref{cap4:con:par}-\eqref{cap4:con:nar} and their linearization ($C=2\sqrt{2}h^2 J$).
			\label{fig:linearization}}
	\end{figure}
	Suppose that 
	constraint $y\le \frac{x}{2} + \frac{C}{2\sqrt{x}}$ is active at $\wb^{(k)}$  (the case when $y\ge \frac{x}{2} - \frac{C}{2\sqrt{x}}$ is active can be dealt with in a completely analogous way). If we linearize such constraint around $\wb^{(k)}$, then we obtain a linear constraint (black line in Figure~\ref{fig:linearization}) which defines a region completely contained into the one defined by the non-linear constraint $y\le \frac{x}{2} + \frac{C}{2\sqrt{x}}$.
	Hence, for each direction~$\dwb^{(k)}$ feasible with respect to the linearized constraint, we are always able to perform sufficiently small steps, without violating the original non-linear constraints, i.e., for $\alpha>0$ small enough,
	it holds that $\wb^{(k)} + \alpha \dwb^{(k)} \in\Omega$.
\end{remark}
A special case occurs when $w_{i-1}+w_{i+1}=0$ holds for some $i\in \{2,\ldots,n-1\}$. In that case, coefficients $\eta_i, b_{P_i}, b_{N_i}$
are not even defined. In fact, in this case we can omit the linearized constraints (\ref{cap4:eq:nar_lin})-(\ref{cap4:eq:par_lin}).
Indeed, the corresponding non-linear constraints \eqref{cap4:con:par}-\eqref{cap4:con:nar}  are not active at the current solution and, thus, along the computed direction a step with strictly positive length can always be taken without violating them.
For all feasible solutions $\wb$ such that this special case does not occur, i.e., such that $w_{i-1}+w_{i+1}>0$, $i=2,\ldots,n-1$, constraints \eqref{cap4:con:par} and \eqref{cap4:con:nar} can be rewritten as follows
\begin{eqnarray}
w_{i-1} + w_{i+1}-2 w_i -2 \sqrt{2} h^2 J (w_{i+1} + w_{i-1})^{-\frac{1}{2}}\leq 0&  \label{cap4:con:par1}\\ [6pt]
2 w_i-w_{i-1} - w_{i+1}-2 \sqrt{2} h^2 J (w_{i+1} + w_{i-1})^{-\frac{1}{2}}\leq 0. & \label{cap4:con:nar1}
\end{eqnarray}
Note that the functions on the left-hand side of these constraints are concave. Now we can define a variant of Problem \ref{cap4:prob_lin} where constraints \eqref{cap4:eq:nar_lin} and~\eqref{cap4:eq:par_lin}  are replaced by
the following linearizations of constraints \eqref{cap4:con:par1} and \eqref{cap4:con:nar1}
\begin{eqnarray}
-\beta_i \dw_{i-1}  - \beta_i \dw_{i+1} +\dw_{i} \le \bni & \label{cap4:con:linnar1} \\ [6pt]
\theta_i \dw_{i-1}  + \theta_i \dw_{i+1} -\dw_{i} \le \bpi, &  \label{cap4:con:linpar1}
\end{eqnarray}
where $\bpi$ and $\bni$ are the same as in  \eqref{cap4:eq:nar_lin} and~\eqref{cap4:eq:par_lin}  (see (\ref{eq:paramdef})), while
\begin{equation}
\label{eq:newcoeff}
\begin{array}{l}
\theta_i=\frac{1}{2}+\frac{\sqrt{2}h^2 J}{2} (w_{i+1} + w_{i-1})^{-\frac{3}{2}} \\ [6pt]
\beta_i=\frac{1}{2}-\frac{\sqrt{2}h^2 J}{2} (w_{i+1} + w_{i-1})^{-\frac{3}{2}}.
\end{array}
\end{equation}
The following proposition states that constraints \eqref{cap4:con:linnar1}-\eqref{cap4:con:linpar1} are tighter than
constraints  \eqref{cap4:eq:nar_lin}-\eqref{cap4:eq:par_lin}.
\begin{prop}
\label{prop:tighter}
For all $i=2,\ldots,n-1$, it holds that $\beta_i\leq \eta_i\leq \theta_i$. Equality $\eta_i=\theta_i$ holds if the corresponding non-linear constraint \eqref{cap4:con:par1} is active at the current point $\wb$. Similarly, $\eta_i=\beta_i$ holds if the corresponding non-linear constraint \eqref{cap4:con:nar1} is active at the current point $\wb$.
\end{prop}
\begin{proof}
We only prove the results about $\theta_i$ and $\eta_i$. Those about $\beta_i$ and $\eta_i$ are proved in a completely analogous way.
By definition of $\eta_i$ and $\theta_i$, we need to prove that
$$
\frac{3(w_{i+1} + w_{i-1}) -2w_i}{4(w_{i+1}+w_{i-1})}\leq \frac{1}{2}+\frac{\sqrt{2}h^2 J}{2} (w_{i+1} + w_{i-1})^{-\frac{3}{2}}.
$$
After few simple computations, this inequality can be rewritten as
$$
w_{i+1} + w_{i-1}-2w_i -2\sqrt{2}h^2 J(w_{i+1} + w_{i-1})^{-\frac{1}{2}}\leq 0,
$$
which holds in view of feasibility of $\wb$ and, moreover, holds as an equality if constraint \eqref{cap4:con:par1} is active at the current point $\wb$, as we wanted to prove.
\end{proof}
In view of this result, by replacing constraints  \eqref{cap4:eq:nar_lin}-\eqref{cap4:eq:par_lin} with \eqref{cap4:con:linnar1}-\eqref{cap4:con:linpar1}, we reduce the search space of the new displacement $\dwb$. On the other hand, the following proposition states that with constraints \eqref{cap4:con:linnar1}-\eqref{cap4:con:linpar1} no line search is needed along
the direction $\dwb$, i.e., we can always choose the step length $\alpha=1$.
\begin{prop}
If constraints \eqref{cap4:con:linnar1}-\eqref{cap4:con:linpar1} are employed as a replacement of
constraints \eqref{cap4:eq:nar_lin}-\eqref{cap4:eq:par_lin} 
in the definition of Problem
\ref{cap4:prob_lin}, then for each feasible solution $\dwb$ of this problem it holds that
$\wb+\dwb\in \Omega$.
\end{prop}
\begin{proof}
For the sake of convenience, let us rewrite Problem \ref{cap4:prob_disc} in the following more compact form
\begin{equation}
\label{eq:origcompact}
\begin{array}{ll}
\min & f({\bf w}+\dwb) \\ [6pt]
& {\bf c}({\bf w}+\dwb)\leq 0,
\end{array}
\end{equation}
where the vector function ${\bf c}$ contains all constraints of Problem \ref{cap4:prob_disc} and
the non-linear ones are given as in \eqref{cap4:con:par1}-\eqref{cap4:con:nar1} (recall that in that case vector ${\bf c}$ is a vector of concave functions).
Then, Problem \ref{cap4:prob_lin} can be written as follows
\begin{equation}
\label{eq:lincompact}
\begin{array}{ll}
\min & f({\bf w}+\delta {\bf w}) \\ [6pt]
& {\bf c}({\bf w}) +\nabla {\bf c}({\bf w}) \delta {\bf w} \leq 0.
\end{array}
\end{equation}
Now, it is enough to observe that, by concavity
$$
{\bf c}({\bf w}+\dwb)\leq {\bf c}({\bf w}) +\nabla {\bf c}({\bf w}) \delta {\bf w},
$$
so that each feasible solution of \eqref{eq:lincompact} is also feasible for \eqref{eq:origcompact}.
\end{proof}
The above proposition states that the feasible region of Problem \ref{cap4:prob_lin}, when constraints
\eqref{cap4:con:linnar1}-\eqref{cap4:con:linpar1} are employed in its definition, is a subset of the feasible region $\Omega$ of
the original Problem \ref{cap4:prob_disc}.
As a final result of this section, we state the following theorem, which establishes convergence of Algorithm SCA
to a stationary (KKT) point of Problem \ref{cap4:prob_disc}, if it runs for an infinite number of iterations and if
constraints
\eqref{cap4:con:linnar1}-\eqref{cap4:con:linpar1} are always employed after a finite number of iterations  in the definition of Problem  \ref{cap4:prob_lin}. 
\begin{thm}
\label{thm:convergence}
If Algorithm SCA is run for an infinite number of iterations and there exists some positive integer value $K$ such that for all iterations $k\geq K$,
constraints
\eqref{cap4:con:linnar1}-\eqref{cap4:con:linpar1} are always employed in the definition of Problem  \ref{cap4:prob_lin},
then the sequence of points $\{{\bf w}^{(k)}\}$ generated by the algorithm converges to a KKT point of Problem  \ref{cap4:prob_disc}.
\end{thm}
\begin{proof}
See Appendix \ref{sec:appconvergence}.
\end{proof}
\begin{remark}
In Algorithm SCA at each iteration we solve to optimality Problem  \ref{cap4:prob_lin}. This is indeed necessary in the final iterations
to prove the convergence result stated in Theorem \ref{thm:convergence}. However, during the first iterations it is not necessary to solve the problem to optimality: finding a feasible descent direction is enough. This does not alter
the theoretical properties of the algorithm and allows to reduce the computing times.
\end{remark}
In the rest of the paper we will refer to constraints~\eqref{cap4:eq:acc_lin}-\eqref{cap4:eq:dec_lin} as acceleration constraints, while constraints \eqref{cap4:eq:nar_lin}-\eqref{cap4:eq:par_lin} (or \eqref{cap4:con:linnar1}-\eqref{cap4:con:linpar1}) will be called (linearized) Negative Acceleration Rate (NAR) and Positive Acceleration Rate (PAR) constraints, respectively.
Also note that in the different subproblems discussed in the following sections we will always refer to the linearization
with constraints \eqref{cap4:eq:nar_lin}-\eqref{cap4:eq:par_lin} and, thus, with parameters $\eta_i$, but the same results also hold for the linearization
with constraints \eqref{cap4:con:linnar1}-\eqref{cap4:con:linpar1} and, thus, with parameters $\theta_i$ and $\beta_i$.
\section{The subproblem with acceleration and NAR constraints}\label{cap4:sec:accnar}
In this section we will propose an efficient method to solve Problem~\ref{cap4:prob_lin} when PAR constraints are removed.
The solution of this subproblem will become part of an approach to solve a suitable relaxation of Problem~\ref{cap4:prob_lin} and, in fact, under very mild assumptions, to solve Problem~\ref{cap4:prob_lin} itself.
This will be clarified in Section \ref{cap4:sec:PAR}.
We will discuss: (i) the subproblem including only (\ref{cap4:eq:bound_lin}) and the acceleration 
constraints~\eqref{cap4:eq:acc_lin} and \eqref{cap4:eq:dec_lin}; (ii) the subproblem including only (\ref{cap4:eq:bound_lin}) and the NAR constraints~\eqref{cap4:eq:nar_lin};
(iii) the subproblem including all constraints~\eqref{cap4:eq:bound_lin}-\eqref{cap4:eq:nar_lin}. 
Throughout the section we will need 
the results stated in the following two propositions.
Let us consider problems with the following form, where $N=\{1,\ldots,n\}$ and $M_j=\{1,\ldots,m_j\}$, $j\in N$,
\begin{equation}
\label{eq:specstruct}
\begin{array}{lll}
\min & g(x_1,\ldots,x_n) & \\ 
 & x_j\leq a_{i,j} x_{j-1} + b_{i,j} x_{j+1} + c_{i,j} &  i\in M_j,\ \ \ j\in N \\ 
& \ell_j\leq x_j \leq u_j & j\in N
\end{array}
\end{equation}
with 
\begin{itemize}
\item $g$ a monotonic decreasing function;
\item $a_{ij}, b_{ij}, c_{ij}\geq 0$, for $ i\in M_j$ and $j\in N$;
\item $a_{i1}=0$ for $i\in M_1$;
\item $b_{in}=0$ for $ i\in M_n$.
\end{itemize}
The following result is proved in~\cite{ConLocLAu19Graph}. Here we report the proof in order to make the paper self-contained.
We denote by $P$ the feasible polytope of problem (\ref{eq:specstruct}) . Moreover, we denote by ${\bf z}$ the component-wise maximum of all feasible solutions in $P$, i.e., for each $j\in N$:
$$
z_j=\max_{{\bf x}\in P} x_j
$$
(note that the above maximum value is attained since $P$ is a polytope).
\begin{prop}
\label{prop:specstruct}
The unique optimal solution  of (\ref{eq:specstruct}) is the component-wise maximum ${\bf z}$ of all its feasible solutions.
\end{prop}
\begin{proof}
If we are able to prove that the component-wise maximum ${\bf z}$ of all feasible solutions is itself a feasible solution, by monotonicity of $g$, it must also be the unique optimal solution. In order to prove that ${\bf z}$ is feasible, 
we proceed as follows. For $j\in N$, let ${\bf x}^{*j}$ be the optimal solution of $\max_{{\bf x}\in P} x_j$, so that
$z_j=x_j^{*j}$. Since ${\bf x}^{*j}\in P$, then it must hold that $\ell_j\leq z_j\leq u_j$. Moreover,
let us consider the generic constraint
$$
 x_j\leq a_{i,j} x_{j-1} + b_{i,j} x_{j+1} + c_{i,j} ,
$$
for $i\in M_j$.
It holds that
$$
\begin{array}{ll}
z_j & =x_j^{*j}\leq  a_{i,j} x_{j-1}^{*j} + b_{i,j} x_{j+1}^{*j} + c_{i,j} \leq \\
 &\leq a_{i,j} z_{j-1}  + b_{i,j} z_{j+1}+ c_{i,j},
\end{array}
$$
where the first inequality follows from feasibility of ${\bf x}^{*j}$, while the second follows from nonnegativity of $a_{ij}$ and $b_{ij}$ and the definition of ${\bf z}$. Since this holds for all $j\in N$, the result is proved.
\end{proof}  
Now, consider the problem obtained from (\ref{eq:specstruct}) by removing some constraints, i.e., by taking $M'_j\subseteq M_j$ for each $j\in N$:
\begin{equation}
\label{eq:specstruct1}
\begin{array}{lll}
\min & g(x_1,\ldots,x_n) & \\ 
 & x_j\leq a_{i,j} x_{j-1} + b_{i,j} x_{j+1} + c_{i,j} &  i\in M'_j,\ \ \ j\in N \\ 
& \ell_j\leq x_j \leq u_j & j\in N,
\end{array}
\end{equation}
Later on we will also need the result stated in the following proposition. 
\begin{prop}
\label{prop:specstruct1}
The optimal solution $\bar{{\bf x}}^\star$ of problem (\ref{eq:specstruct1}) is an upper bound for the optimal solution 
${\bf x}^\star$ of problem (\ref{eq:specstruct}), i.e., $\bar{{\bf x}}^\star\geq {\bf x}^\star$.
\end{prop}
\begin{proof}
It holds that ${\bf x}^\star$ is a feasible solution of problem (\ref{eq:specstruct1}), so that, in view of Proposition
\ref{prop:specstruct}, $\bar{{\bf x}}^\star\geq {\bf x}^\star$ holds.
\end{proof}
\subsection{Acceleration constraints}
The simplest case is the one where we only consider the acceleration constraints ~\eqref{cap4:eq:acc_lin} and~\eqref{cap4:eq:dec_lin}, besides constraints~\eqref{cap4:eq:bound_lin} with a generic upper bound vector ${\bf y}\ge {\bf 0}$. The problem to be solved is:
\begin{problem}
	\label{cap4:prob_lin_acc}
	\begin{align*}
	\min_{\dwb \in \Real^n}& \sum_{i=1}^{n-1} \frac{2h}{\sqrt{w_{i+1} +\dw_{i+1}} + \sqrt{w_{i}+\dw_{i}}}\\
	& \mathbf{l_B}  \leq  \dwb  \leq \mathbf{y}, \\
	& \dw_{i+1}-\dw_{i}  \le \bAi, & i=1,\dots,n-1,\\
	&\dw_{i}-\dw_{i+1}  \le \bddi, & i=1,\dots n-1.
	\end{align*}
\end{problem}
It can be seen that such problem belongs to the class of problems (\ref{eq:specstruct}).
Therefore, in view of Proposition \ref{prop:specstruct}, the optimal solution of Problem \ref{cap4:prob_lin_acc} is the component-wise maximum of its feasible region.
Moreover, in ~\cite{consolini2017scl} it has been proved that
 Algorithm~\ref{cap4:alg:acc}, based on a forward and a backward iteration and with $O(n)$ computational complexity, returns an optimal solution of Problem \ref{cap4:prob_lin_acc}.
\begin{algorithm}
	\caption{Routine \texttt{SolveAcc} for the solution of the problem with acceleration constraints\label{cap4:alg:acc}}
	\SetKwData{Left}{left}\SetKwData{This}{this}\SetKwData{Up}{up}
	\SetKwFunction{Union}{Union}\SetKwFunction{FindCompress}{FindCompress}
	\SetKwInOut{Input}{input}\SetKwInOut{Output}{output}
	\Input{Upper bound $\mathbf{y}$ }
	\Output{$\dwb$}
	$\dw_1 = 0$, $\dw_n = 0$ \;
	
	\For{$i=1$ \KwTo $n-1$}{
		
		$\dw_{i+1} = \min \begin{Bmatrix}{\dw_{i} + \bAi}, y_{i+1}\end{Bmatrix}$
	}	
	\For{$i=n-1$ \KwTo $1$}{
		
		$\dw_{i} = \min \begin{Bmatrix}{\dw_{i+1} +\bAi}, y_{i}\end{Bmatrix}$
	}
	\Return $\dwb$
\end{algorithm}
\subsection{ NAR constraints }
Now, we consider the problem only including NAR constraints~\eqref{cap4:eq:nar_lin} and constraints~\eqref{cap4:eq:bound_lin} with upper bound vector ${\bf y}$:
\begin{problem}\label{cap4:prob:nar}
	\begin{align}
	& \qquad \min_{\dwb \in \Real^n} \sum_{i=1}^{n-1} \frac{2h}{\sqrt{w_{i+1} +\dw_{i+1}} + \sqrt{w_{i}+\dw_{i}}} \nonumber
	\\
	&  \mathbf{0}  \leq  \dwb  \leq \mathbf{y},&  \label{cap4:con:nar_bound} \\
	&\dw_{i} \le \eta_i (\dw_{i-1} + \dw_{i+1}) + \bni , &i = 2,\dots,n-1,\label{cap4:con:nar_nar}
	\end{align}
\end{problem}
where $y_1=y_n=0$ because of the boundary conditions.
Also this problem belongs to the class of problems (\ref{eq:specstruct}), so that Proposition \ref{prop:specstruct} states that
its optimal solution is the component-wise maximum of its feasible region.
Problem \ref{cap4:prob:nar} can be solved by using the graph-based approach presented in~\cite{CabConLoc2018coap1,ConLocLAu19Graph}.
However, reference~\cite{CabConLoc2018coap1} shows that, by exploiting the structure of a simpler version of the NAR constraints, it is possible to develop an algorithm more efficient than the graph-based one.
Our purpose is to extend the results presented in reference~\cite{CabConLoc2018coap1} to a case with different and more challenging NAR constraints, in order to develop an efficient algorithm outperforming
the graph-based one. 
%
\newline\newline\noindent

Now, let us consider the restriction of Problem \ref{cap4:prob:nar} between two generic indexes $s,t$ such that $1\leq s<t\leq n$, obtained by fixing $\dw_s = y_s$ and $\dw_t = y_t$ and by considering only the NAR and upper bound constraints at $s+1,\ldots,t-1$.
Let $\dwb^*$ be the optimal solution of the restriction.
We first prove the following lemma.
\begin{lem}\label{cap4:rem:max}
The optimal solution $\dwb^*$ of the restriction of Problem \ref{cap4:prob:nar} between two indexes $s,t$, $1\leq s<t\leq n$, is such that
for each $j\in \{s+1,\ldots,t-1\}$, either $\dw_j^* \le y_j$ or
	$\dw_j^* \le \eta_j(\dw^*_{j+1} + \dw^*_{j-1}) + b_{N_j}$ holds as an equality.
\end{lem}
\begin{proof}
It is enough to observe that in case both inequalities were strict for some $j$, then, in view of the monotonicity of the objective function, we could decrease the objective function value by increasing the value of $\dw^*_j$, thus contradicting optimality of $\dwb^*$.
\end{proof}
Note that  the above result also applies to the full Problem \ref{cap4:prob:nar}, which corresponds to the special case $s=1$, $t=n$ with $y_1=y_n=0$.
In view of Lemma~\ref{cap4:rem:max} we have that there exists an index $j$, with  $s< j \le t$, such that: (i) $\dw^*_j = y_j$;  
(ii) the upper bound constraint is not active at $s+1,\ldots,j-1$; (iii) all NAR constraints $s+1,\dots,j-1$ are active.
Then, $j$ is the lowest index in $\{s+1,\ldots, t-1\}$ where the upper bound constraint is active
If index $j$ were known, then the following observation allows to return the components of the optimal solution between $s$ and $j$.
Let us first introduce the following definitions of matrix ${\bf A}$ and vector ${\bf q}$:
\begin{equation}
\label{eq:defA}	
	{\bf A}= \begin{bmatrix}
	1 & -\eta_{s+1} &  0   & \cdots & 0 \\ 
	-\eta_{s+2} & 1 &-\eta_{s+2}& \ddots& \vdots\\ 
	0 & \ddots & \ddots& \ddots & 0  \\
	\vdots& \ddots & -\eta_{j-2} & 1 & -\eta_{j-2} \\
	0 & \cdots& 0  & -\eta_{j-1} & 1  
	\end{bmatrix}\,, 
\end{equation}
\begin{equation}
\label{eq:defq}
{\bf q}
	= \begin{bmatrix}{\bn}_{s+1} + \eta_{s+1}y_s \\{\bn}_{s+2}\\ \vdots\\
	{\bn}_{j-2}\\ {\bn}_{j-1} + \eta_{j-1}y_j\end{bmatrix}.
	\end{equation}
	Note that {\bf A} is the square submatrix of the NAR constraints restricted to rows $s+1$ up to $j-1$ and the related columns.
\begin{obser}
	Let $\dwb^*$ be the optimal solution of the restriction of Problem \ref{cap4:prob:nar} between $s$ and $t$ and let $s < j$.
	If constraints $\dw^*_s\le y_s$, $\dw^*_j\le y_j$, and $\dw^*_i \le \eta_i(\dw^*_{i+1} + \dw^*_{i-1}) + b_{N_i}$, for $i = s+1,\dots,j-1$, are all active, then $\dw^*_{s+1},\dots,\dw^*_{j-1}$ are obtained by the solution of the following
tridiagonal system
$$
\begin{array}{ll}
\dw_s=y_s & \\
\dw_{r}-\eta_{r}\dw_{r+1}-\eta_{r}\dw_{r-1}= {\bn}_{r}  & r=s+1,\ldots,j-1 \\ 
\dw_j=y_j, &
\end{array}
$$
or, equivalently, as 
\begin{equation}
\label{eq:linsysaux}
\begin{array}{ll}
\dw_{s+1}-\eta_{s+1}\bar{x}_{s+2}= {\bn}_{s+1} +\eta_{s+1} y_s & \\ 
\dw_{r}-\eta_{r}\dw_{r+1}-\eta_{r}\dw_{r-1}= {\bn}_{r}  & r=s+2,\ldots,j-2 \\ 
\dw_{s+1}-\eta_{s+1}\bar{x}_{s+2}= {\bn}_{s+1} +\eta_{s+1} y_s. 
\end{array}
\end{equation}
In matrix form, the above tridiagonal linear system can be written as follows, where matrix ${\bf A}$ is defined in (\ref{eq:defA}), while vector ${\bf q}$ is defined in (\ref{eq:defq}):	
	\begin{equation}\label{linear-system}
	{\bf A} \begin{bmatrix} \dw^*_{s+1}\\ \vdots \\ \dw^*_{j-1} \end{bmatrix}   = {\bf q}.
	\end{equation}
\end{obser}
Tridiagonal systems 
\[	a_i x_{i-1} + b_i x_i + c_i x_{i+1} = d_i, \ \ \ i=1,\ldots,m,	\]
with $a_1  = c_m = 0$, can be solved through so called Thomas algorithm~\cite{higham2002accuracy} (see Algorithm~\ref{cap4:algo:thomas}) with  $O(m)$ operations.
\begin{algorithm}[!h]
	\caption{Thomas algorithm}\label{cap4:algo:thomas}
	\SetKwInOut{Input}{input}\SetKwInOut{Output}{output}
	\Input{$\mathbf{a}$, $\mathbf{b}$, $\mathbf{c}$, $\mathbf{d}$}
	\Output{$\bar\xb$}
	Let $m$ be the dimension of ${\bf d}$\\
	\tcc{ {\em Forward phase}}\
	\For{$i=2,\dots,m$}
	{
		Set $\delta_i = \frac{a_i}{b_{i-1}}$\;
		Set $b_i = b_i - \delta_i c_{i-1}$\;
		Set $d_i = d_i - \delta_id_{i-1}$\;
		Set $\alpha_i = \frac{d_i}{b_i} $\;\label{cap4:algo:thomas:alpha}
		Set $\psi_i = \frac{c_i}{b_i}$\;\label{cap4:algo:thomas:beta}
	}
	\tcc{ {\em Backward phase}}\
	Set $\bar{x}_m = \alpha_m$\; \label{cap4:alg:thomas:back_start}
	\For {$i = m-1,\dots,1$} 
	{
		Set $\bar{x}_i = \alpha_i - \psi_i \bar{x}_{i+1}$ \label{cap4:alg:thomas:back_end}\;
	}
\end{algorithm}
In order to detect the lowest index $j\in \{s+1,\ldots, t-1\}$ such that the upper bound constraint is active at $j$, 
we propose Algorithm~\ref{cap4:algo:solveNAR}, also called \texttt{SolveNAR}  and described in what follows.
We initially set $j=t$. Then, at each iteration we 
solve the linear system~\eqref{linear-system}.
Let $\bar{\bf x} = (\bar{x}_{s+1},\dots,\bar{x}_{j-1})$ be its solution.
We check whether it is feasible  and optimal or not. 
Namely, if there exists $k\in \{s+1,\ldots,j-1\}$ such that either $\bar{x}_k< 0$ or $\bar{x}_k > y_k$, then $\bar{\bf x}$ is unfeasible and, consequently, we need to reduce $j$ by 1.
%
%
%
If $\bar{x}_k= y_k$ for some $k\in\{s+1,\ldots,j-1\}$, then we also reduce $j$ by 1 since $j$ is not in any case the lowest index of the optimal solution where the upper bound constraint is active.
Finally, if  $0\le \bar{x}_k< y_k$, for $k= s+1,\dots,j-1$, 
then we need to verify if 
${\bf \bar{x}}$ is the best possible
solution over the interval $\{s+1,\dots,j-1\}$. We will be able to check that after
proving the following result.
\begin{prop}
\label{prop:optimalpiece}
Let matrix ${\bf A}$ be defined as in (\ref{eq:defA}) and vector ${\bf q}$ be defined as in (\ref{eq:defq}).
The optimal solution $\dwb^*$ of the restriction of Problem \ref{cap4:prob:nar} between $s$ and $t$  satisfies
\begin{equation}
\label{eq:pieceopt}
\dw_s^*=y_s,\ \ \dw_r^*=\bar{x}_r,\ \ r=s+1,\ldots,j-1,\ \ \dw_j^*=y_j, 
\end{equation}
if and only if the optimal value of the LP problem:
\begin{equation}\label{cap4:prob:epsilon_red}
\begin{array}{cl}
\max_{\boldsymbol{\epsilon}} & \boldsymbol{1}^T\boldsymbol{\epsilon}\\
&{\bf A} \boldsymbol{\epsilon} \le {\bf 0}, \\
& \boldsymbol{\epsilon} \leq {\bf \bar y} - {\bf\bar \xb},
\end{array}
\end{equation}
is strictly positive or, equivalently, if the following system admits no solution:
\begin{equation}\label{cap4:KKT}
{\bf A}^T \boldsymbol{\lambda} = {\bf 1}, \  \boldsymbol{\lambda} \ge {\bf 0}.
\end{equation}
\end{prop}
\begin{proof}
Let us first assume that  $\dwb^*$ does not fulfill (\ref{eq:pieceopt}). Then, in view of Lemma \ref{cap4:rem:max},
$j$ is not the lowest index such that the upper bound is active at the optimal solution and, consequently, $\dw^*_k=y_k>\bar{x}_k$ for some $k\in \{s+1,\ldots,j-1\}$. 
Such optimal solution must be feasible and, in particular, it must satisfy
all NAR constraints  between $s+1$ and $j-1$ and the upper bound constraints between $s+1$ and $j$, i.e.:
$$
\begin{array}{ll}
\dw^*_{s+1}-\eta_{s+1}\dw^*_{s+2}\leq {\bn}_{s+1} +\eta_{s+1} y_s & \\ 
\dw^*_{r}-\eta_{r}\dw^*_{r+1}-\eta_{r}\dw^*_{r-1}\leq {\bn}_{r}  & r=s+2,\ldots,j-2 \\ 
\dw^*_{j-1}-\eta_{j-1}\dw^*_{j-2}-\eta_{j-1}\dw^*_j \leq {\bn}_{j-1} & \\ 
\dw^*_{r}\leq y_r & r=s+1,\ldots,j. 
\end{array}
$$
In view of $\dw^*_j\leq y_j$ and $\eta_{j-1}\geq 0$, $\dwb^*$ also satisfies the following system of inequalities:
$$
\begin{array}{ll}
\dw^*_{s+1}-\eta_{s+1}\dw^*_{s+2}\leq {\bn}_{s+1} +\eta_{s+1} y_s & \\ 
\dw^*_{r}-\eta_{r}\dw^*_{r+1}-\eta_{r}\dw^*_{r-1}\leq {\bn}_{r}  & r=s+2,\ldots,j-2 \\ 
\dw^*_{j-1}-\eta_{j-1}\dw^*_{j-2} \leq {\bn}_{j-1} +\eta_{j-1} y_j & \\ 
\dw^*_{r}\leq y_r & r=s+1,\ldots,j-1. 
\end{array}
$$
After making the change of variables $\dw^*_r=\bar{x}_r+\epsilon_r$ for $r=s+1,\ldots,j-1$, and
recalling that ${\bar {\bf x}}$ solves system (\ref{eq:linsysaux}),
the system of inequalities can be further rewritten as:
$$
\begin{array}{ll}
\epsilon_{s+1}-\eta_{s+1}\epsilon_{s+2}\leq 0 & \\ 
\epsilon_{r}-\eta_{r}\epsilon_{r+1}-\eta_{r}\epsilon_{r-1}\leq 0 & r=s+2,\ldots,j-2 \\ 
\epsilon_{j-1}-\eta_{j-1}\epsilon_{j-2} \leq 0 & \\ 
\epsilon_{r}\leq y_r - \bar{x}_{r}& r=s+1,\ldots,j-1. 
\end{array}
$$
Finally, recalling the definition of matrix ${\bf A}$ and vector ${\bf q}$ given in (\ref{eq:defA}) and (\ref{eq:defq}), respectively, this can also be written in a more compact form as:
$$
\begin{array}{l}
{\bf A}\boldsymbol{\epsilon}\leq {\bf 0} \\ 
\boldsymbol{\epsilon}\leq \bar{{\bf y}}-\bar{{\bf x}}.
\end{array}
$$
If $\dw^*_k=y_k>\bar{x}_k$ for some $k\in \{s+1,\ldots,j-1\}$, then the system must admit a solution with $\epsilon_k>0$. This is equivalent to prove that problem (\ref{cap4:prob:epsilon_red})
has an optimal solution with at least one strictly positive component and the optimal value is strictly positive. Indeed, in view of the definition of matrix ${\bf A}$, problem (\ref{cap4:prob:epsilon_red}) has the structure of the problems discussed in Proposition \ref{prop:specstruct}. More precisely, to see that we need to remark that maximizing
$\boldsymbol{1}^T\boldsymbol{\epsilon}$ is equivalent to minimizing the decreasing function $-\boldsymbol{1}^T\boldsymbol{\epsilon}$. Then, observing that $\boldsymbol{\epsilon}={\bf 0}$ is a feasible solution of problem (\ref{cap4:prob:epsilon_red}), by Proposition \ref{prop:specstruct} the optimal solution $\boldsymbol{\epsilon}^*$ must be a nonnegative vector, and since at least one component, namely component $k$, is strictly positive, then the optimal value must also be strictly positive.
\newline\newline\noindent
Conversely, let us assume that the optimal value is strictly positive and that $\boldsymbol{\epsilon}^*$ is an optimal solution with at least one strictly positive component. Then, there are two possible alternatives. Either
the optimal solution  $\dwb^*$ of
the restriction of Problem \ref{cap4:prob:nar} between $s$ and $t$ is such that $\dw_j^*<y_j$, in which case (\ref{eq:pieceopt}) obviously does not hold,
or $\dw_j^*=y_j$. In the latter case, let us assume by contradiction that (\ref{eq:pieceopt}) holds. We observe that the solution defined as follows:
$$
\begin{array}{ll}
x'_s=y_s & \\
x'_r=\bar{x}_r+\epsilon^*_r=\dw_r^*+ \epsilon^*_r& r=s+1,\ldots, j-1 \\
x'_j=y_j =\dw_j^*& \\
x'_r=\dw_r^* & r=j+1,\ldots,t,
\end{array}
$$
is feasible for the restriction of Problem \ref{cap4:prob:nar} between $s$ and $t$. Indeed, by feasibility of $\boldsymbol{\epsilon}^*$ in problem (\ref{cap4:prob:epsilon_red}) all upper bound and NAR constraints between $s$ and $j-1$ are fulfilled. Those between, $j+1$ and $t$ are also fulfilled by the feasibility of $\dwb^*$. Then, we only need to prove that the NAR constraint at $j$ is satisfied. By feasibility of $\dwb^*$ and in view of $\epsilon_{j-1}^*,\eta_j\geq 0$, we have that :
$$
\begin{array}{l}
x'_j=\dw_j^*\leq \eta_{j} \dw_{j-1}^*+\eta_{j} \dw_{j+1}^*+{\bn}_{j}\leq \\
\leq  \eta_{j} (\dw_{j-1}^*+\epsilon_{j-1})+\eta_{j} \dw_{j+1}^*+{\bn}_{j}=\eta_{j} x'_{j-1}+\eta_{j} x'_{j+1}+{\bn}_{j}.
\end{array}
$$
Thus, ${\bf x}'$ is feasible and such that ${\bf x}'\geq \dwb^*$ with at least one strict inequality (recall that at least one component of $\boldsymbol{\epsilon}^*$ is strictly positive), which contradicts the optimality of $\dwb^*$ (recall that
the optimal solution must be the component-wise maximum of all feasible solutions).
\newline\newline\noindent
In order to prove the last part, i.e., that problem  (\ref{cap4:prob:epsilon_red})  has a positive optimal value if and only if (\ref{cap4:KKT}) admits no solution, we notice that the optimal value is positive if and only if the feasible point ${ \boldsymbol{\epsilon} = {\bf 0}}$ is not an optimal solution, or, equivalently, the null vector is not a KKT point.
Since at  ${\boldsymbol \epsilon = {\bf 0}}$  constraints $\boldsymbol{\epsilon} \le {\bf \bar y} - {\bf\bar \xb}$ cannot be active, then the KKT conditions for problem~\eqref{cap4:prob:epsilon_red} at this point are exactly those established in
(\ref{cap4:KKT}), where vector $\boldsymbol{\lambda}$ s the vector of Lagrange mutlpliers for constraints ${\bf A} \boldsymbol{\epsilon} \le {\bf 0}$. This concludes the proof.
\end{proof}
Then if~\eqref{cap4:KKT} admits no solution, then
(\ref{eq:pieceopt}) does not hold and, again, we need 
to reduce $j$ by 1. Otherwise, we can fix the optimal solution between $s$ and $j$ according to (\ref{eq:pieceopt}). 
After that, we recursively call the routine \texttt{SolveNAR} on the remaining subinterval $\{j,\dots,t\}$ in order to obtain the solution
over the full interval.
\begin{remark}
	In Algorithm~\ref{cap4:algo:solveNAR} routine \texttt{isFeasible} is the routine used to verify if, for $ k = s+1,\dots,j-1$, $0\le \bar{x}_k < y_k$, while \texttt{isOptimal} is the procedure to check optimality of $\bar{ \bf x}$ over the interval $\{s+1,\dots,j-1\}$, i.e., that (\ref{eq:pieceopt}) holds.
\end{remark}
Now, we are ready to prove that Algorithm~\ref{cap4:algo:solveNAR} solves Problem \ref{cap4:prob:nar}.
\begin{prop}
	The call \texttt{solveNAR($\by$,1,$n$)} of  Algorithm~\ref{cap4:algo:solveNAR}  returns the optimal solution of
	Problem \ref{cap4:prob:nar}.
\end{prop}
\begin{proof}
After the call \texttt{solveNAR($\by$,1,$n$)}, we are able to identify the portion of the optimal solution between 1 and some index $j_1$, $1<j_1\leq n$.
If $j_1=n$, then we are done. Otherwise, we make the recursive call \texttt{solveNAR($\by$,$j_1$,$n$)}, which will enable
 to identify also the portion of the optimal solution between $j_1$ and some index $j_2$, $j_1<j_2\leq n$.
If $j_2=n$, then we are done. Otherwise, we make the recursive call \texttt{solveNAR($\by$,$j_2$,$n$)}, and so on. After at most $n$ recursive calls, we are able to return the full optimal solution.
\end{proof}
\begin{algorithm}[!ht]
	\caption{\texttt{SolveNAR}(${\bf y}, s, t$)}\label{cap4:algo:solveNAR}
	\SetKwFunction{Union}{Union}\SetKwFunction{FindCompress}{FindCompress}
	\SetKwInOut{Input}{input}\SetKwInOut{Output}{output}
	\Input{Upper bound $\mathbf{y}$ and two indices $s$ and  $t$ with $1\le s<t \le n$}
	\Output{$\dwb^*$ }
	Set $j = t$\;
	$\dwb^* = {\bf y}$\;
	
	\While{$j \ge s+1$}{
		
		Compute the solution $\bar{\bf x}$ of the linear system~\eqref{linear-system}\;
		\If{
			\texttt{isFeasible}$(\bar{\bf x})$ and \texttt{isOptimal}$(\bar{\bf x})$}{
			\textbf{Break}\;
		}
		\Else{
			Set $j = j- 1$\;
			
		}
	}
	\For{$i = s+1,\dots,j-1$}{
		Set $\dw^*_i = \bar{x}_{i}$\;}
	
	\Return $\dwb^* = \min \left\{\dwb^*, \texttt{SolveNAR}(\dwb^*,j,t)  \right\} $\;
\end{algorithm}
\begin{remark}
Note that Algorithm~\ref{cap4:algo:solveNAR} involves solving
a significant amount of linear systems, both to compute $\bar{\bf x}$ and to verify its optimality (see~\eqref{linear-system} and~\eqref{cap4:KKT}).
In what follows we propose some implementation details which improve the performance of Algorithm~\ref{cap4:algo:solveNAR}.
As previously remarked, each (tridiagonal) linear system~\eqref{linear-system} can be solved in at most $O(j-s)$ operations (here the number $m$ of equations is $t-s-1$), which can be upper bounded by $O(n)$.
Moreover, we can directly check the feasibility of $\bar{\bf x}$ during the backward phase of the Thomas Algorithm (see lines~\ref{cap4:alg:thomas:back_start}-\ref{cap4:alg:thomas:back_end} of Algorithm~\ref{cap4:algo:thomas}), namely we declare unfeasibility as soon as $0\le \bar{x}_i\le y _i$ does not hold, without completing the backward propagation. We also observe that coefficients
$\alpha_i$ and $\psi_i$, $i=2,\ldots,m$ do not change with $j$, so that the forward phase of the Thomas algorithm can be performed only once at the beginning of the procedure \texttt{solveNAR} for the whole interval $\{s,\dots,t\}$.
Finally, Thomas algorithm can also be employed to solve the (tridiagonal) linear system~\eqref{cap4:KKT}, needed to verify optimality of $\bar{\bf x}$.
It is also worthwhile to remark that $j$ can be reduced by more than one unit at each iteration. Indeed, let $m_{i,r}=\prod_{s=i}^{r-1} \psi_s$ for $i<r\leq j$. Then, it can be seen
that $\bar{x}_i=q_r-m_{i,r} \bar{x}_r$, for some $q_r$ and each $r\in \{i+1,\ldots,j\}$. Now, let us assume that $\bar{x}_i> y_i$. In such case
what we are currently doing is moving form $j$ to $j-1$ and compute a new solution $\bar{{\bf x}}^{{\tt new}}$. However, we are able to compute in advance the value
 $\bar{x}_i^{{\tt new}}$ without solving the full linear system. Indeed, we have
$$
\bar{x}_i^{{\tt new}}=\bar{x}_i+m_{i,j-1} (\bar{x}_{j-1}-y_{j-1}),
$$
and, in case $\bar{x}_i^{{\tt new}}> y_i$, we can further reduce to $j-2$  and repeat the same procedure. A similar approach can be employed when $\bar{x}_i<0$.
\end{remark}
The following proposition states the worst-case complexity of \texttt{solveNAR($\by$,1,$n$)}.
\begin{prop}
	Problem \ref{cap4:prob:nar} can be solved with $O(n^3)$ operations by running the procedure \texttt{SolveNAR}$({\bf y},1,n)$  and by using the Thomas algorithm for the solution of each linear system.
\end{prop}
\begin{proof}
	In the worst case, at the first call we have $j_1=2$, since we need to go all the way from $j=n$ down to $j=2$.
Since for each $j$ we need to solve a tridiagonal system, which requires at most $O(n)$ operations, the first call of \texttt{SolveNAR} requires $O(n^2)$ operations. This is similar for all successive calls, and since the number of recursive calls is at most $O(n)$, the overall effort is at most of $O(n^3)$ operations.
\end{proof}
In fact, what we observed is that the practical complexity of the algorithm is much better, namely $\Theta(n^2)$.
\subsection{Acceleration and NAR constraints}
Now we discuss the problem with acceleration and NAR constraints, with upper bound vector ${\bf y}$, i.e.:
\begin{problem}
	\label{cap4:prob_lin_accnar}
	\begin{align*}
	& \qquad \min_{\dwb \in \Real^n} \sum_{i=1}^{n-1} \frac{2h}{\sqrt{w_{i+1} +\dw_{i+1}} + \sqrt{w_{i}+\dw_{i}}}
	\\
	& \mathbf{l_B}  \leq  \dwb  \leq \mathbf{y}, \\
	& \dw_{i+1}-\dw_{i}  \le \bAi, & i=1,\dots,n-1,\\
	&\dw_{i}-\dw_{i+1}  \le \bddi, & i=1,\dots n-1,\\
	&\dw_{i} - \eta_i \dw_{i-1}  - \eta_i \dw_{i+1} \le \bni ,& i = 2,\dots,n-1.
	\end{align*}
\end{problem} 
We first remark that Problem \ref{cap4:prob_lin_accnar} has the structure of problem (\ref{eq:specstruct}), so that by Proposition
\ref{prop:specstruct}, its unique optimal solution is the component-wise maximum of its
	feasible region.
As for Problem \ref{cap4:prob:nar}, we can solve  Problem \ref{cap4:prob_lin_accnar}
by using the graph-based approach proposed in reference~\cite{ConLocLAu19Graph}.
However, reference~\cite{CabConLoc2018coap1} shows that, if we adopt a very efficient procedure
to solve Problems \ref{cap4:prob_lin_acc} and  \ref{cap4:prob:nar}, 
then it is worth splitting
the full problem into two separated ones and 
use an iterative approach (see Algorithm~\ref{cap4:alg:accnar}). 
Indeed, Problems \ref{cap4:prob_lin_acc}-\ref{cap4:prob_lin_accnar} share the common property that their optimal solution is also the component-wise maximum of the corresponding feasible region. Moreover,
according to Proposition \ref{prop:specstruct1}, the optimal solutions of 
Problems \ref{cap4:prob_lin_acc} and  \ref{cap4:prob:nar} are valid upper bounds for the optimal solution (actually, also for any feasible solution) of the full Problem \ref{cap4:prob_lin_accnar}.
In Algorithm~\ref{cap4:alg:accnar} we first call the procedure $\texttt{SolveACC}$ with input the upper bound
vector ${\bf y}$. Then, the output of this procedure, which, according to what we have just stated, is an upper bound
 for the solution of the full Problem \ref{cap4:prob_lin_accnar}, satisfies $\dwb_{\text{Acc}}\leq {\bf y}$ and becomes the input for a call of the procedure
$ \texttt{SolveNAR}$. The output $\dwb_{\text{NAR}}$ of this call will be again an upper bound for the solution of the full Problem \ref{cap4:prob_lin_accnar} and it will satisfy $\dwb_{\text{NAR}}\leq \dwb_{\text{ACC}}$. This output will become the input of a further call to the procedure $\texttt{SolveACC}$, and we proceed in this way until the distance between two consecutive output vectors falls below a prescribed tolerance value $\varepsilon$. The following proposition states that the sequence of output vectors generated by the alternate calls to the procedures $\texttt{SolveACC}$ and $ \texttt{SolveNAR}$ will converge to the optimal solution of the full Problem \ref{cap4:prob_lin_accnar}.
\begin{prop}
	Algorithm~\ref{cap4:alg:accnar} converges to the the optimal solution of Problem \ref{cap4:prob_lin_accnar} when $\varepsilon=0$ and stops after a finite number of iterations if $\varepsilon>0$.
\end{prop}
\begin{proof}
We have observed that the sequence of alternate solutions of Problems \ref{cap4:prob_lin_acc} and  \ref{cap4:prob:nar}, here denoted by $\{{\bf y}_t\}$, is:
(i) a sequence of valid upper bounds for the optimal solution of Problem \ref{cap4:prob_lin_accnar}; (ii) component-wise monotonic non-increasing; (iii) component-wise bounded from below by the null vector. 
Thus, if $\varepsilon=0$, an infinite sequence is generated which
converges to some point $\bar{{\bf y}}$, which is also an upper bound for the optimal solution of Problem \ref{cap4:prob_lin_accnar} but, more precisely, by continuity is also a feasible point of the problem and, is thus, also the optimal solution
of the problem. If $\varepsilon>0$, due to the convergence to some point $\bar{{\bf y}}$, at some finite iteration the exit condition of the while loop must be satisfied.
\end{proof}
\begin{algorithm}[h]
	\caption{ Algorithm \texttt{SolveACCNAR} for the solution of  Problem \ref{cap4:prob_lin_accnar}.}
	\label{cap4:alg:accnar}
	\SetKwData{Left}{left}\SetKwData{This}{this}\SetKwData{Up}{up}
	\SetKwFunction{Union}{Union}\SetKwFunction{FindCompress}{FindCompress}
	\SetKwInOut{Input}{input}\SetKwInOut{Output}{output}
	\Input{ The upper bound $\mathbf{y}$ and the tolerance $\varepsilon$}
	\Output{The optimal solution $\dwb^*$ and the optimal value $f^{*}$}
	$\dwb_{\text{Acc}} = $ \texttt{SolveACC}$(\mathbf{y})$\;
      $\dwb_{\text{NAR}}= $ \texttt{SolveNAR}$(\dwb_{\text{Acc}},1,n)$\;
	\While{$\|\dwb_{\text{NAR}}  - \dwb_{\text{Acc}} \| > \varepsilon$}{
		
		$\dwb_{\text{Acc}} = $  \texttt{SolveACC}$(\dwb^*)$\;
            $\dwb_{\text{NAR}} = $ \texttt{SolveNAR}$(\dwb_{\text{Acc}},1,n)$\;
	}	
$\dwb^*=\dwb_{\text{NAR}}$\;
\Return $\dwb^*$, \texttt{evaluateObj}$(\dwb^*)$
\end{algorithm}

\section{A descent method for the case of acceleration, PAR and NAR constraints}\label{cap4:sec:PAR}
\begin{figure*}[!h]
	\centering
     \includegraphics[width=\textwidth]{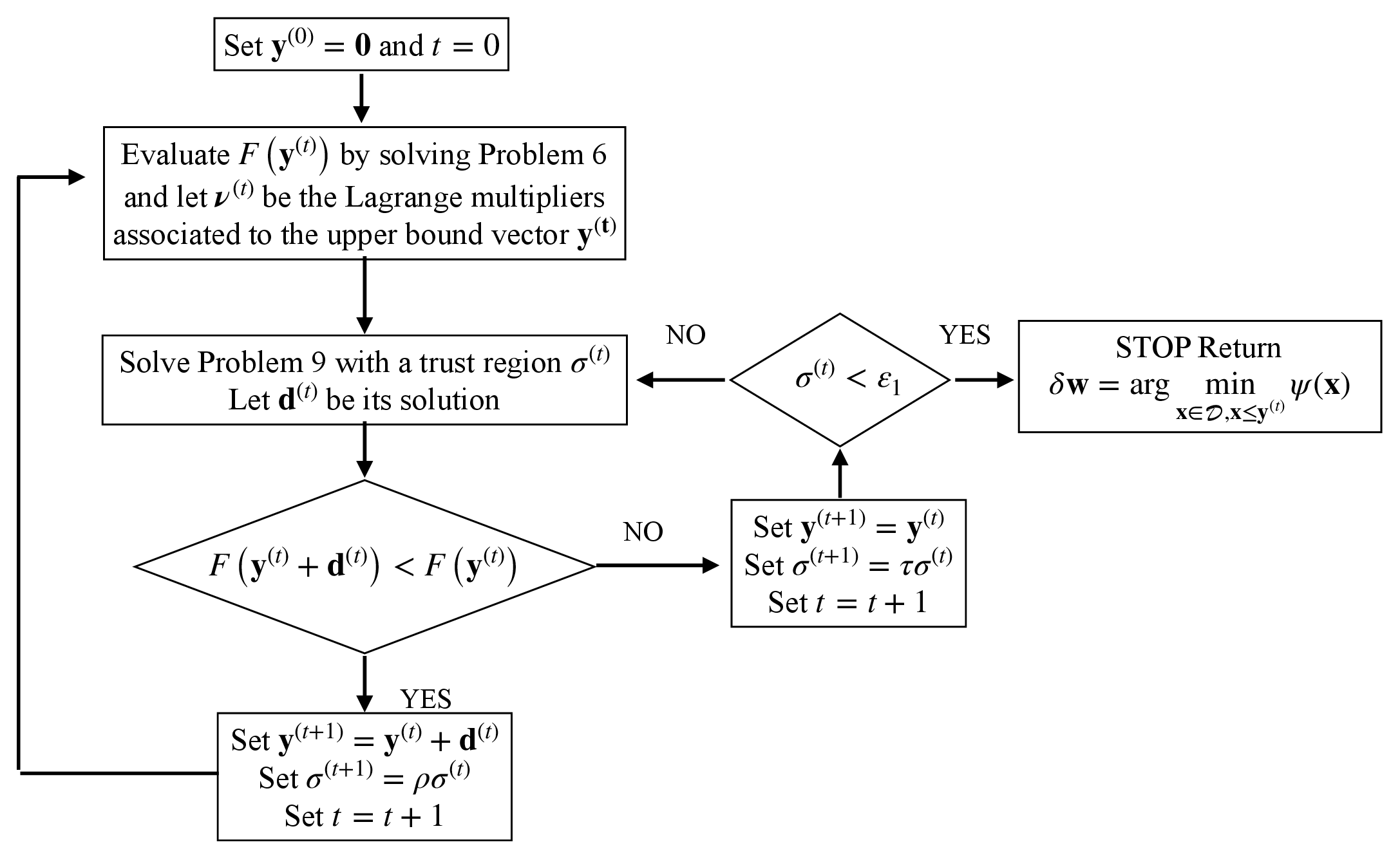}
	\caption{Flow chart of the routine \texttt{ComputeUpdate}
		\label{fig:flowcompute}}
\end{figure*}

Unfortunately, PAR constraints~\eqref{cap4:eq:par_lin}
do not satisfy the assumptions requested in Proposition
\ref{prop:specstruct} in order to guarantee that the component-wise maximum of the feasible region is the optimal solution of Problem~\ref{cap4:prob_lin}.
However, in Section~\ref{cap4:sec:accnar} we have shown that  Problem \ref{cap4:prob_lin_accnar}, i.e., Problem~\ref{cap4:prob_lin}  without the PAR constraints, can be efficiently solved by Algorithm~\ref{cap4:alg:accnar}.
Our purpose then is to separate
the acceleration and NAR constraints from the  PAR constraints.
\begin{defn}\label{cap4:defn:optimvalfun}
	Let $f: \Real^n\rightarrow\Real$ be the objective function of Problem~\ref{cap4:prob_lin} and let $\mathcal{D}$  be the region defined by the acceleration and NAR constraints (the feasible region of Problem \ref{cap4:prob_lin_accnar}). We define the function $F: \Real^n \rightarrow \Real$ as follows
	\[
	F(\yb) = \min\left\{ f(\xb)\, | \,  \xb \in \mathcal{D}, \xb \le \yb \right\}.
	\]
	Namely, $F$ is the optimal value function of Problem \ref{cap4:prob_lin_accnar}
	when the upper bound vector is $\yb$.
\end{defn}
\begin{prop}
	Function $F$ is a convex function.
\end{prop}
\begin{proof}
	Since  Problem \ref{cap4:prob_lin_accnar} is convex, then the optimal
	value function $F$ is convex (see Section 5.6.1 of~\cite{boyd2004convex}) .
\end{proof}
Now, let us introduce the following problem:
\begin{problem}\label{cap4:prob:par}
	\begin{align}
	&\quad \min_{\by \in \Real^n} F(\by) \label{cap4:con:par_obj}\\
	\intertext{such that}
	& \eta_i( y_{i-1}  +  y_{i+1}) -y_{i} \le {\bpi} ,\quad i = 2,\dots,n-1, \label{cap4:con:par_par}\\
	& \mathbf{l_B}\le \by \le \mathbf{u_B} \label{cap4:con:par_bound}.
	\end{align}
\end{problem}
Such problem is a relaxation of Problem~\ref{cap4:prob_lin}. Indeed, each feasible solution of Problem~\ref{cap4:prob_lin} is also feasible for Problem~\ref{cap4:prob:par} and the value of $F$ at such solution is equal to the value of the objective function of Problem~\ref{cap4:prob_lin} at the same solution. 
We will solve Problem~\ref{cap4:prob:par} rather than  Problem~\ref{cap4:prob_lin} to compute the new displacement  $\dwb$. More precisely, if $\yb^*$ is the optimal solution of Problem~\ref{cap4:prob:par}, then we will set 
\begin{equation}
\label{eq:computedir}
\dwb  = \arg \min_{\xb \in \mathcal{D}, \xb \le \by^*} f(\xb).
\end{equation}
In the following proposition we prove that, under a very mild condition, the optimal solution of Problem \ref{cap4:prob:par}
computed in (\ref{eq:computedir}) is feasible and, thus, optimal for Problem \ref{cap4:prob_lin}, so that,
although we solve a relaxation of the latter problem, we return an optimal solution for it.
\begin{prop}\label{cap4:prop:salvezza}
Let $\wb^{(k)}$ be the current point. 
If
\begin{equation}\label{cap4:ass:prop_salvezza}
\dw_{j-1} + \dw_{j+1}\le 2\left(w^{(k)}_{j-1} + w^{(k)}_{j+1} \right), \quad j=2,\ldots,n-1,
\end{equation}
where $\dwb$ is computed through (\ref{eq:computedir}), 
then, $\dwb$ is feasible for Problem \ref{cap4:prob_lin}, both if the non-linear constraints are linearized as
in  \eqref{cap4:eq:nar_lin}-\eqref{cap4:eq:par_lin} and if they are linearized as in \eqref{cap4:con:linnar1}-\eqref{cap4:con:linpar1}.
\end{prop}
\begin{proof}
See Appendix \ref{sec:feasdir}. 
\end{proof}
Note that assumption~\eqref{cap4:ass:prop_salvezza} is mild since we are basically requiring that no steps larger than twice as much as the current values can be taken. 
In order to fulfill it, one can 
impose restrictions on $\dw_{j+1}$ and $\dw_{j-1}$, like, e.g., 
$$\dw_{j-1} \le w^{(k)}_{j-1} + \frac{w^{(k)}_{j-1} +w^{(k)}_{j+1} }{2},$$
and a similar restriction for $\dw_{j+1}$, so that the
assumption is satisfied. In fact, in the computational experiments
we did not impose such restrictions unless a positive step-length along the computed direction $\dwb$ could not be taken (which, however, never occurred in our experiments).
\newline\newline\noindent
Now, let us turn our attention towards the solution of Problem~\ref{cap4:prob:par}. In order to solve it, we propose a descent method. 
We can exploit the information provided by
the dual optimal solution $\nub\in\Real_{+}^n$ associated to the upper bound constraints of Problem \ref{cap4:prob_lin_accnar}.
Indeed, from the sensitivity theory, we know that the dual solution is related
to the gradient of the optimal value function $F$ (see Definition~\ref{cap4:defn:optimvalfun}) and provides information about
how it changes its value for small perturbations of the upper bound values (for further details see Sections 5.6.2 and 5.6.5 in~\cite{boyd2004convex}). 
%
Let $\by^{(t)}$ be  a feasible solution of Problem~\ref{cap4:prob:par} and $\nub\in\Real_{+}^n$ be the Lagrange multipliers of the upper bound constraints of  Problem \ref{cap4:prob_lin_accnar}, when
the upper bound is $\yb^{(t)}$.
Let:
$$\varphi_i = \bpi - \eta_i\left(  y^{(t)}_{i-1}  + y^{(t)}_{i+1} \right) + y^{(t)}_i, \quad i=2,\dots,n-1.$$
Then, a \emph{feasible descent direction} $\bd^{(t)}$ 
can be obtained by solving the following LP problem: 
\begin{problem}\label{cap4:prob:direction}
	\begin{align}
	&\min_{\bd \in \Real^n} -\nub^{T}\bd \label{cap4:con:obj_dir}\\
	& \eta_i\left(  d_{i-1}  + d_{i+1} \right)- d_i \le \varphi_i, \, i = 2,\dots,n-1, \label{cap4:con:par_dir}\\
	& \mathbf{l_B} \le \yb^{(t)} + \bd \le \mathbf{u_{B}} \label{cap4:con:box_dir},
	\end{align}
\end{problem}
where the objective function~\eqref{cap4:con:obj_dir} imposes that $\bd^{(t)}$ is a descent direction while constraints~\eqref{cap4:con:par_dir} and~\eqref{cap4:con:box_dir} guarantee feasibility with respect to Problem~\ref{cap4:prob:par}. 
Problem~\ref{cap4:prob:direction} is an LP problem and, consequently, it can easily be solved through a standard LP solver. In particular, we employed GUROBI \cite{gurobi}. 
Unfortunately, since the information provided by the dual optimal solution $\nub$ is local and related to
small perturbations of the upper bounds, it might happen that $F(\by^{(t)} + \bd^{(t)}) \ge F(\by^{(t)})$. 
To overcome this issue we introduce a trust-region constraint in Problem~\ref{cap4:prob:direction}. So,
let $\sigma^{(t)}\in\Real_{+}$ be  the radius of the trust-region at iteration $t$, then we have:
\begin{problem}\label{cap4:prob:direction_tr}
	\begin{align}
	&\min_{\bd \in \Real^n} -\nub^{T}\bd \label{cap4:obj:d}\\
	& \eta_i\left(  d_{i-1}  + d_{i+1} \right)- d_i \le \varphi_i, \,  i = 2,\dots,n-1, \label{cap4:con:par_dir_tr}\\
	&\mathbf{\bar{l}_{B}} \le \bd \le \mathbf{\bar{u}_B} \label{cap4:con:box_dir_tr}.
	\end{align}
\end{problem}
where $\bar{l}_{B_i} = \max\{l_{B_i}-y^{(t)}_i,-\sigma^{(t)}\}$ and $\bar{u}_{B_i} = \min\{u_{B_i}-y_i^{(t)},\sigma^{(t)}\}$ for $i=1,\dots,n$.
After each iteration of the descent algorithm, we change the radius $\sigma^{(t)}$ 
according to the following rules:
\begin{itemize}
	\item if $F(\by^{(t)}+\bd^{(t)}) \ge F(\by^{(t)})$,
	then we set $\by^{(t+1)} = \by^{(t)}$ and we tight the trust-region by decreasing $\sigma^{(t)}$ by a factor $\tau\in(0,1)$;
	\item if~$F(\by^{(t)}+\bd^{(t)}) < F(\by^{(t)})$,
	then we set $\by^{(t+1)} = \by^{(t)} + \bd^{(t)}$ and enlarge the radius~$\sigma^{(t)}$ by a factor $\rho>1$. 
\end{itemize}
The proposed descent algorithm is sketched in Figure~\ref{fig:flowcompute}, which reports the flow chart of the procedure~\texttt{ComputeUpdate} used in Algorithm SCA. 
We initially set $\by^{(0)} = {\bf 0}$. At each iteration $t$ we evaluate the objective function $F(\by^{t})$ by solving 
Problem \ref{cap4:prob_lin_accnar} with upper bound vector $\by^{(t)}$
through a call of the routine \texttt{solveACCNAR} (Algorithm~\ref{cap4:alg:accnar}). Then, we compute the Lagrange multipliers
$\nub^{(t)}$ associated to the upper bound constraints.
After that, we compute a candidate descent direction $\bd^{(t)}$ by solving Problem~\ref{cap4:prob:direction_tr}.
If $\bd^{(t)}$ is a descent step, then we set $\by^{(t+1)} = \by^{(t)} + \bd^{(t)}$  and enlarge the radius of the trust region, otherwise we do not move to a new point and we tight the trust region and solve again~Problem~\ref{cap4:prob:direction_tr}.
The descent algorithm stops as soon as the radius of the trust region becomes smaller than a fixed tolerance $\varepsilon_1$.
\begin{remark}
	Note that we initially set $\by^{(0)} ={\bf 0} $. But any feasible solution of Problem~\ref{cap4:prob:direction_tr} does the job and, actually, starting with a good initial solution may enhance the performance of the algorithm.
\end{remark}
\begin{remark}
\label{rem:heuristic}
Problem~\ref{cap4:prob:direction_tr} is an LP one and can be solved by any existing LP solver. However,
a 
suboptimal solution to Problem~\ref{cap4:prob:direction_tr}, obtained by a heuristic approach, is also acceptable. Indeed, we observe that:
i) an {\em optimal} descent direction is not strictly required; ii) a heuristic approach allows to reduce the time needed to get a descent direction.
In this paper we propose a possible heuristic. This will be described in Appendix~\ref{sec:heuristic}. However, we point out that a possible topic for future research is the development of further heuristic approaches.
\end{remark}

\section{Speed planning in general configuration spaces}
\label{sec_higher_dim}
In this section, we consider the speed planning problem for a curve in
a generic configuration space. We show that, also in this case, a speed profile obtained by
solving Problem~\ref{cap4:prob:continuos} allows to bound the
velocity, the acceleration and the jerk of the obtained trajectory.
Let $\mathcal{Q}$ be a smooth manifold of dimension $p$ that represents a
configuration space with $p$-degrees of
freedom ($p$-DOF). 
For instance, the configuration space of a rigid body
corresponds to $SE(3)$, the set of rigid transformations in $\Real^3$.
Let $\|\cdot\|: T Q \to \Real$ be a a norm on $T \mathcal{Q}$, the
tangent space of $\mathcal{Q}$.
Let $\bgamma: [0,s_f] \to \mathcal{Q}$ be a $C^3$ function,
whose image set $\Gamma=\bgamma([0,s_f])$ is the path to be followed
and such that $\bgamma$ has unit-length parameterization,
that is $(\forall s \in [0,s_f]) \|\bgamma'(s)\|=1$. In this way,
$s_f$ is the length of $\Gamma$.
In particular, $\bgamma(0)$ and $\bgamma(s_f)$ are the initial and
final configurations. Define $t_{f}$ as the time when the
configuration reaches the end of the path. Let
 $\lambda : [0, t_f] \rightarrow [0, s_f]$ be a differentiable
 monotone increasing function such that $\gamma(\lambda(t))$ is the
 configuration at time $t$ and let   $ v : [0, s_f] \rightarrow [0, +\infty]$ be such that
  $\left( \forall t \in [0,t_f]\right) \dot{\lambda}(t) =
  v(\lambda(t))$. Namely, $v(s)$ is the norm of the velocity of the
 configuration along $\Gamma$ at position $\bgamma(s)$. We impose ($\forall s \in [0,s_{f}]$) $v(s) \ge 0$.
For  any $t \in [0,t_f]$, using chain rule, we obtain, setting $w=v^2$
and $\bq(t) = \bgamma(\lambda(t))$,
\begin{equation}  
\label{eq:rep}
\begin{array}{ll}
\dot{\bq}(t) =  \bgamma^{\prime}(\lambda(t))v(\lambda(t)),\\[3pt]
\ddot{\bq}(t) =
                                                                     \frac{1}{2} \bgamma^{\prime}(\lambda(t))w^\prime(\lambda(t))+ \bgamma^{\prime\prime}(\lambda(t))w(\lambda(t)),\\
  \dddot{\bq}(t) =\frac{3}{2}\bgamma^{\prime
  \prime}(\lambda(t))w^\prime(\lambda(t))w(\lambda(t))^\frac{1}{2}+\\
  \bgamma^{\prime}(\lambda(t))w^{\prime \prime}(\lambda(t)) w(\lambda(t))^\frac{1}{2}+ 
  \bgamma^{\prime\prime \prime}(\lambda(t))w(\lambda(t))^\frac{3}{2}\,.
\end{array}
\end{equation}
At this point, one could formulate a speed optimization problem,
similar in structure to Problem~\ref{cap4:prob:continuos}, with constraints
on $\|\dot{\bq}(t) \|$, $\|\ddot{\bq}(t)\|$, $\|\dddot{\bq}(t)\|$.
This leads to a different optimization problem that, although related
to Problem~\ref{cap4:prob:continuos}, would need a different
discussion and is outside the scope of this paper.
However, in the following, we show that if the speed profile $w$ is
chosen by solving Problem~\ref{cap4:prob:continuos}, then
quantities $\|\dot{\bq}(t) \|$, $\|\ddot{\bq}(t)\|$,
$\|\dddot{\bq}(t)\|$ are bounded by terms depending on the parameters
$\mu^+$, $A$, $J$ appearing in Problem~\ref{cap4:prob:continuos}.
To this end, set
$k(\lambda)=\|\bgamma^{\prime\prime}(\lambda(t))\|$,
$k_2(\lambda)=\|\bgamma^{\prime\prime\prime}(\lambda(t))\|$, and, recalling that $\|\gamma'(\lambda(t))\|=1$, note
that 
\[
\begin{array}{ll}
\|\dot{\bq}(t) \| =  v(\lambda(t))\\[3pt]
\|\ddot{\bq}(t)\| \leq \frac{1}{2} |w^\prime(\lambda(t))|+
k(\lambda(t)) w(\lambda(t)),\\
\|\dddot{\bq}(t)\| \leq \frac{3}{2} k(\lambda(t)) |w^\prime(\lambda(t))|  w(\lambda(t))^\frac{1}{2}+
  |w^{\prime \prime}(\lambda(t))| w(\lambda(t))^\frac{1}{2}\\+ 
  k_2(\lambda(t))w(\lambda(t))^\frac{3}{2}\,.
\end{array}
\]

Hence, if $w$ is feasible for
Problem~\ref{cap4:prob:continuos}, the following bounds hold
\[
  \begin{array}{ll}
\|\dot{\bq}(t) \| \leq  \sqrt{\mu^+(\lambda(t))}\\[3pt]
\|\ddot{\bq}(t)\| \leq A+
k(\lambda(t)) \mu^+(\lambda(t)),\\
\|\dddot{\bq}(t)\| \leq 3  k(\lambda(t)) A \mu^+(\lambda(t))^\frac{1}{2}+
  2 J + 
  k_2(\lambda(t))\mu^+(\lambda(t))^\frac{3}{2}\,.
\end{array}
    \]

Hence, a speed profile $w$ obtained as the solution of
Problem~\ref{cap4:prob:continuos} implies that quantities
$\|\dot{\bq}(t) \|$, $\|\ddot{\bq}(t) \|$, $\|\dddot{\bq}(t) \|$ are
bounded in a known way.
If one wants to satisfy constraints
$\|\dot{\bq}(t) \|\leq \hat V$, $\|\ddot{\bq}(t)\| \leq \hat A$,
$\|\dddot{\bq}(t)\| \leq \hat J$, it is possible to proceed in the
following way. Set two constant $0<A<\hat A$, $0<J<J^+$ (for instance,
set $A=\frac{\hat A}{2}$ and $J=\frac{\hat J}{2}$ and define
$\mu^+(\lambda) = \min\{\hat{V}^2, \frac{\hat A-A}{
  k(\lambda)},\chi(\lambda)\}$,
where $\chi(\lambda)$ is a positive quantity that satisfies equation $3  k(\lambda(t)) A \chi(\lambda) ^\frac{1}{2}+
  J+ 
  k_2(\lambda(t)) \chi(\lambda) ^\frac{3}{2} = \hat  J$, then any $w$
  obtained by solving Problem~\ref{cap4:prob:continuos} satisfies the
  required bounds.

\section{Computational Experiments}
\label{sec:compexp}
In this section we present various computational experiments performed in order to evaluate the  approaches proposed in Sections~\ref{cap4:sec:accnar} and~\ref{cap4:sec:PAR}. 
\newline\newline\noindent
In particular, we compared solutions of Problem~\ref{cap4:prob_disc}
computed by algorithm SCA to solutions
obtained with commercial NLP solvers. Note that, with a single exception, we did not
carry out a direct comparison with other methods specifically tailored
to Problem~\ref{cap4:prob_disc} for the following reasons.
\begin{itemize}
\item
Some algorithms (such as
\cite{Villagra-et-al2012}, \cite{RaiCGL2019jerk}) use heuristics to quickly find suboptimal solutions
of acceptable quality, but do not achieve local optimality. Hence
comparing their solution times with SCA would not be fair. However, in one of our experiments (Experiment 4), we made a comparison between the most recent heuristic proposed in
\cite{RaiCGL2019jerk} and Algorithm SCA, both in terms of computing times and in terms of the quality of the returned solution.
\item The method presented in~\cite{8569414}  does not consider the (nonconvex)
  jerk constraint, but solves a convex problem whose objective
  function has a penalization term
  that includes pseudo-jerk. Due to this difference, a direct comparison with SCA is not
  possible.
\item The method presented in~\cite{DONG20071941} 
is based on the numerical solution of a large number of non-linear and non-convex
subproblems and is therefore structurally slower than SCA, whose main
iteration is based on the efficient solution of the convex Problem~\ref{cap4:prob_lin}.
\end{itemize}

In the first two experiments we compare the computational time of IPOPT, a general purpose NLP solver~\cite{ipopt2006}, with that of Algorithm SCA
over some randomly generated instances of Problem~\ref{cap4:prob_disc}.
In particular, we tested two different versions of Algorithm SCA. 
The first version, called \emph{SCA-H} in what follows, employs the heuristic  mentioned in Remark \ref{rem:heuristic} and described
in Appendix~\ref{sec:heuristic}. Since the heuristic procedure may fail in some cases,  in such cases we also need  an LP solver. In particular, in our experiments, we 
used GUROBI whenever the heuristic did not produce either a feasible solution to Problem~\ref{cap4:prob:direction_tr} or a descent direction.
In the second version, called~\emph{SCA-G} in what follows, we always employed GUROBI to solve Problem~\ref{cap4:prob:direction_tr}.
For what concerns the choice of the NLP solver IPOPT, we remark that we chose it after a comparison with two further general purpose NLP solvers, SNOPT and MINOS, which, however, turned out to perform worse than IPOPT on this class of problems.
\newline\newline\noindent
In the third experiment we compare the performance of the two implemented versions of Algorithm SCA 
applied to two specific paths and see their behaviour as the number $n$ of discretized points increases.
\newline\newline\noindent
In the fourth experiment, we compare the solutions returned by Algorithm SCA with those returned by the heuristic recently proposed in \cite{RaiCGL2019jerk}.
\newline\newline\noindent
Finally, in the fifth experiment, in order to illustrate the approach presented in Section~\ref{sec_higher_dim},
we consider a speed planning problem for a UAV vehicle. 
\newline\newline\noindent
We remark that we have also made some experiments to compare the computational time of routine~\texttt{solveACCNAR} (Algorithm~\ref{cap4:alg:accnar}) with the graph-based approach
proposed in~\cite{ConLocLAu19Graph} and with GUROBI for solving Problem \ref{cap4:prob_lin_accnar}. Note that, strictly speaking, Problem \ref{cap4:prob_lin_accnar} is not an LP one since its objective function
is not linear. However, as discussed in~\cite{ConLocLAu19Graph}, its (monotonic non-increasing) objective function can be converted into a (monotonic non-increasing) linear function, thus making GUROBI a valid option to solve the problem.
The computational experiments show that routine~\texttt{solveACCNAR} strongly outperforms both the graph-based approach and GUROBI. That was expected, since 
the graph-based approach and GUROBI are general purpose for a wide class of problems, while routine~\texttt{solveACCNAR} is tailored to the problem with acceleration and NAR constraints.
\newline\newline\noindent
Finally, we remark that rather than employing an NLP solver only once to solve the non-convex Problem~\ref{cap4:prob_disc}, we could have employed it to solve the convex Problem~\ref{cap4:prob_lin} arising at each iteration of the proposed method in place of the procedure \texttt{ComputeUpdate}, presented in this paper. However, the experiments revealed that in doing this the computing times become much larger even with respect to the single call to the NLP solver for solving  the non-convex Problem~\ref{cap4:prob_disc}. This confirms that 
the problem-specific procedure \texttt{ComputeUpdate} is able to strongly outperform a general-purpose NLP solver when solving the convex Problem~\ref{cap4:prob_lin}.
\newline\newline\noindent
All tests have been performed on an  IntelCore i7-8550U CPU at 1.8 GHz. 
Both for IPOPT and Algorithm SCA 
the null vector was chosen as a starting point. The parameters used within Algorithm SCA
were $\varepsilon = 1e^{-8}$, $\varepsilon_1 = 1e^{-6}$ (tolerance parameters), $\rho = 4$ and $\tau = 0.25$ (trust-region update parameters). The initial trust region radius $\sigma^{(0)}$ was initialized to 1 in the first iteration $k=0$, but
adaptively set equal to the size of the last update $\|\wb^{(k)}-\wb^{(k-1)}\|_{\infty}$ in all subsequent iterations (this adaptive choice allowed to reduce computing times by more than a half).
We remark that Algorithm SCA
has been implemented in MATLAB, so we expect better performance after a C/C++ implementation.
\paragraph{Experiment 1}
As a first experiment we compared  the performance of Algorithm SCA
with the NLP solver IPOPT.
We made the experiments over a set of 50 different paths, each of which
was discretized setting $n = 100$, $n = 500$ and $n = 1000$ sample points.
The instances were generated by assuming that the traversed path was divided into
seven intervals over which the curvature of the path was assumed to be constant.
Thus, the $n$-dimensional upper bound vector ${\bf u}$ was generated as follows.
First, we fixed $u_1 = u_n = 0$, i.e., the initial and final speed must be equal to 0.
Next, we partitioned the set $\{2,\dots,n-1\}$ into seven subintervals $I_j, \ j \in \{1,\ldots,7\}$, which
correspond to intervals with constant curvature.
Then, for each subinterval we randomly generated a value $u_j\in(0,\tilde{u}]$, where $\tilde{u}$  is the maximum upper bound (which was set equal to 100 m$^2$s$^{-2}$). Finally, for each $j\in \{1,...,7\}$ we set
$u_k = \tilde{u}_j $ $  \forall k \in I_{j} $.
The  maximum acceleration parameter $A$ is set equal to  $2.78$ ms$^{-2}$
and the maximum jerk $J$ to  0.5
ms$^{-3}$, while the path length is $s_f$ = 60 m.
The values for $A$ and $J$ allow a comfortable motion
for a ground transportation vehicle (see~\cite{hoberock1977survey}).

\begin{figure}
	\centering
	\begin{subfigure}[b]{\columnwidth}
		\includegraphics[width=0.6\columnwidth]{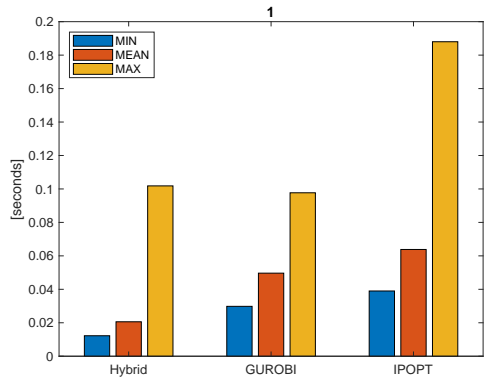}
		\caption{Samples $n = 100$}
		\label{fig:gull}
	\end{subfigure}
	~ 
\centering
	\begin{subfigure}[b]{\columnwidth}
		\includegraphics[width=0.6\columnwidth]{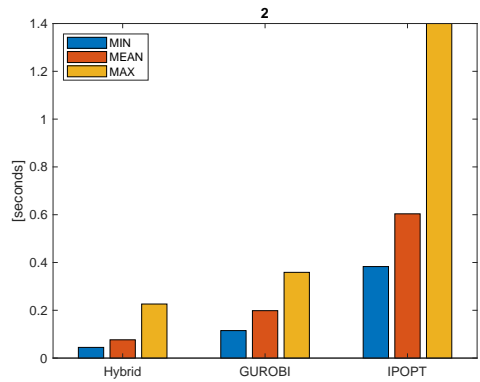}
		\caption{Samples $n = 500$}
		\label{fig:tiger}
	\end{subfigure}
	~ 
\centering
	\begin{subfigure}[b]{\columnwidth}
		\includegraphics[width=0.6\columnwidth]{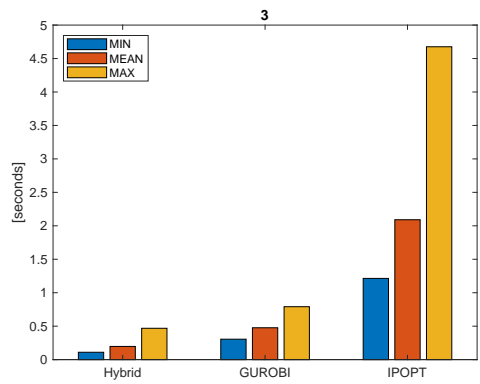}
		\caption{Samples $n = 1000$}
		\label{fig:mouse}
	\end{subfigure}
	\caption{Computational results for Experiment 1.\label{cap:fig:ex2_result}}
\end{figure}

The results are reported in Figure~\ref{cap:fig:ex2_result} in which we
show the minimum, maximum and  average computational times.
The results show that Algorithm \emph{SCA-H} is the fastest one, while \emph{SCA-G} is slightly faster than IPOPT at $n=100$ but clearly faster for a larger number of sample points $n$. In general, we observe that both \emph{SCA-H} and
\emph{SCA-G} tend to outperform IPOPT as $n$ increases.
For what concerns the objective function values returned by the three algorithms, there are some differences due to numerical issues related to the choice of the tolerance parameters, but such differences are mild ones and  never exceed 1\%. 

\begin{figure}[!h]
	\centering
	\includegraphics[width=\columnwidth]{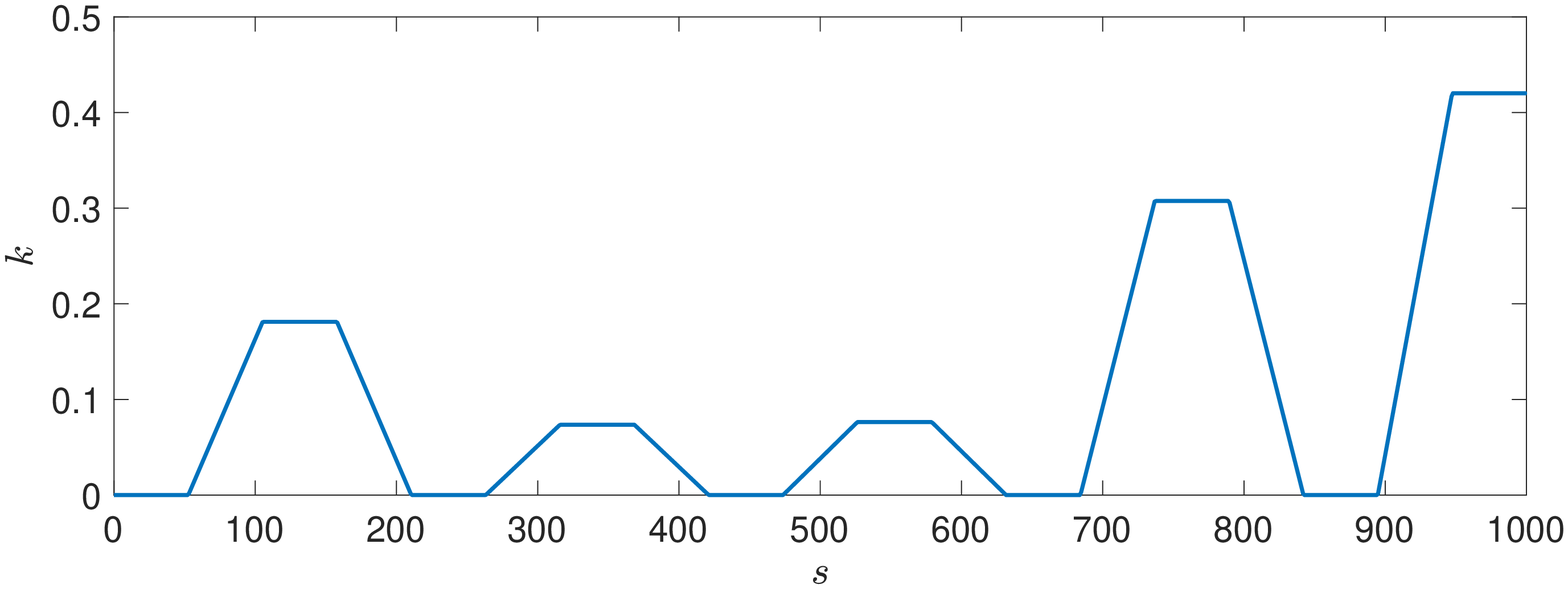}
	\caption{Example of a randomly generated curvature.\label{cap4:fig:curvature}}
	\label{fig:curvatureexample}
\end{figure}

\paragraph{Experiment 2}
We compared again the performance of Algorithm SCA
with the NLP solver IPOPT over different paths.
Again, we made the experiments over a set of 50 different paths, each of which
was discretized using $n = 100$, $n = 500$ and $n = 1000$ variables.
These new instances were randomly generated such that the traversed path was divided into up to five intervals over which the curvature could be zero, linear with respect to the arc-length or constant. We chose this kind of paths since they are able to represent the curvature of a road trip (see \cite{FraSch:04}). An example of the generated curvature is shown in Figure~\ref{cap4:fig:curvature}.
The maximum squared velocity along the path was fixed equal to 192.93 m$^2$s$^{-2}$ (corresponding to a maximum velocity of 50kmh$^{-1}$).
The total length of the paths was fixed to $s_f = 1000$ m, while parameter $A$ was set equal to 0.25 ms$^{-2}$, $J$ to 0.025 ms$^{-3}$ and $A_N$ to 4.9 ms$^{-2}$.
The results are reported in Figure~\ref{cap:fig:ex2_result} in which we
display the minimum, maximum and average computational times.
\begin{figure}
	\centering
	\begin{subfigure}[b]{\columnwidth}
		\includegraphics[width=0.6\columnwidth]{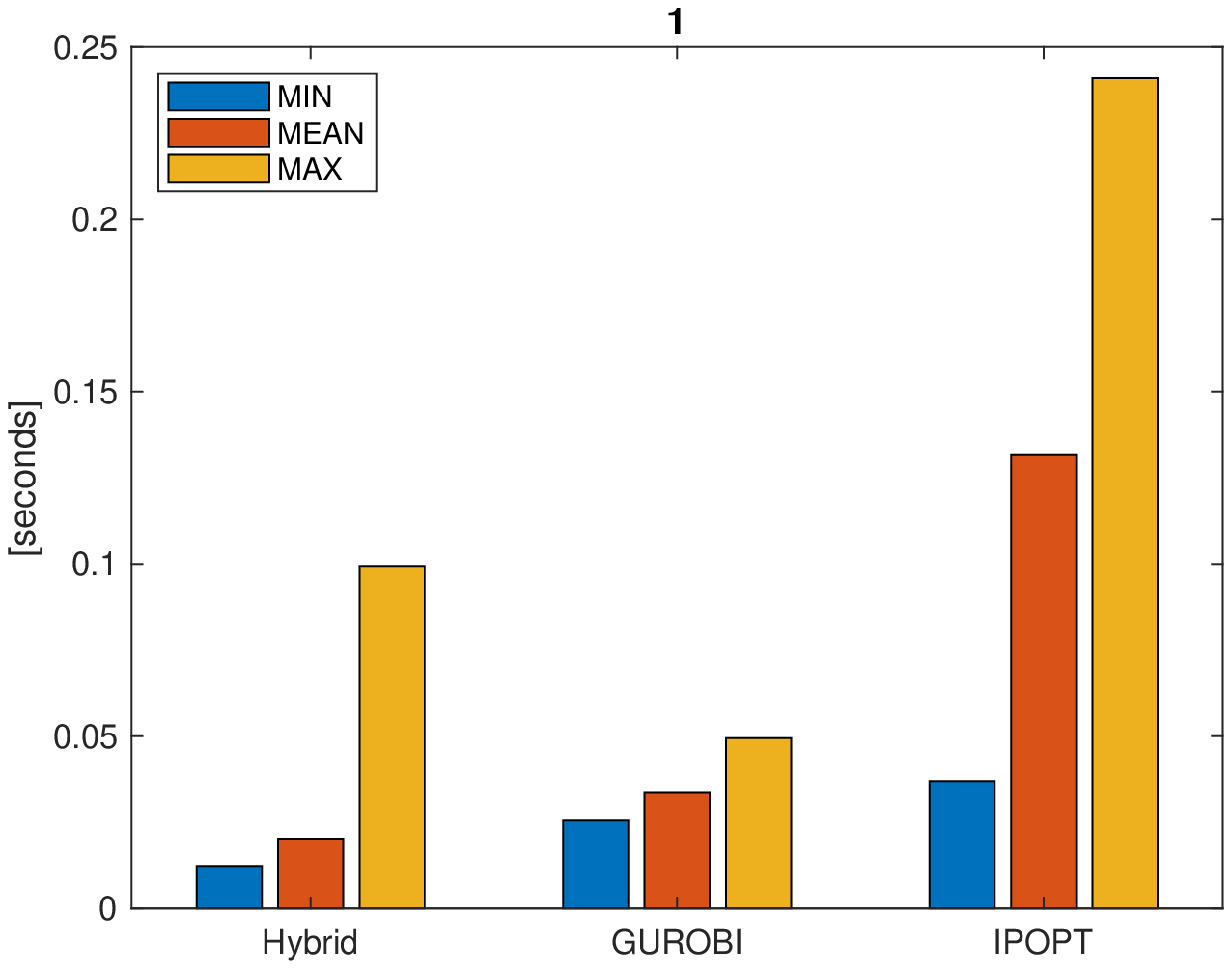}
		\caption{Samples $n = 100$}
		\label{fig:gull}
	\end{subfigure}
	~ 
	\begin{subfigure}[b]{\columnwidth}
		\includegraphics[width=0.6\columnwidth]{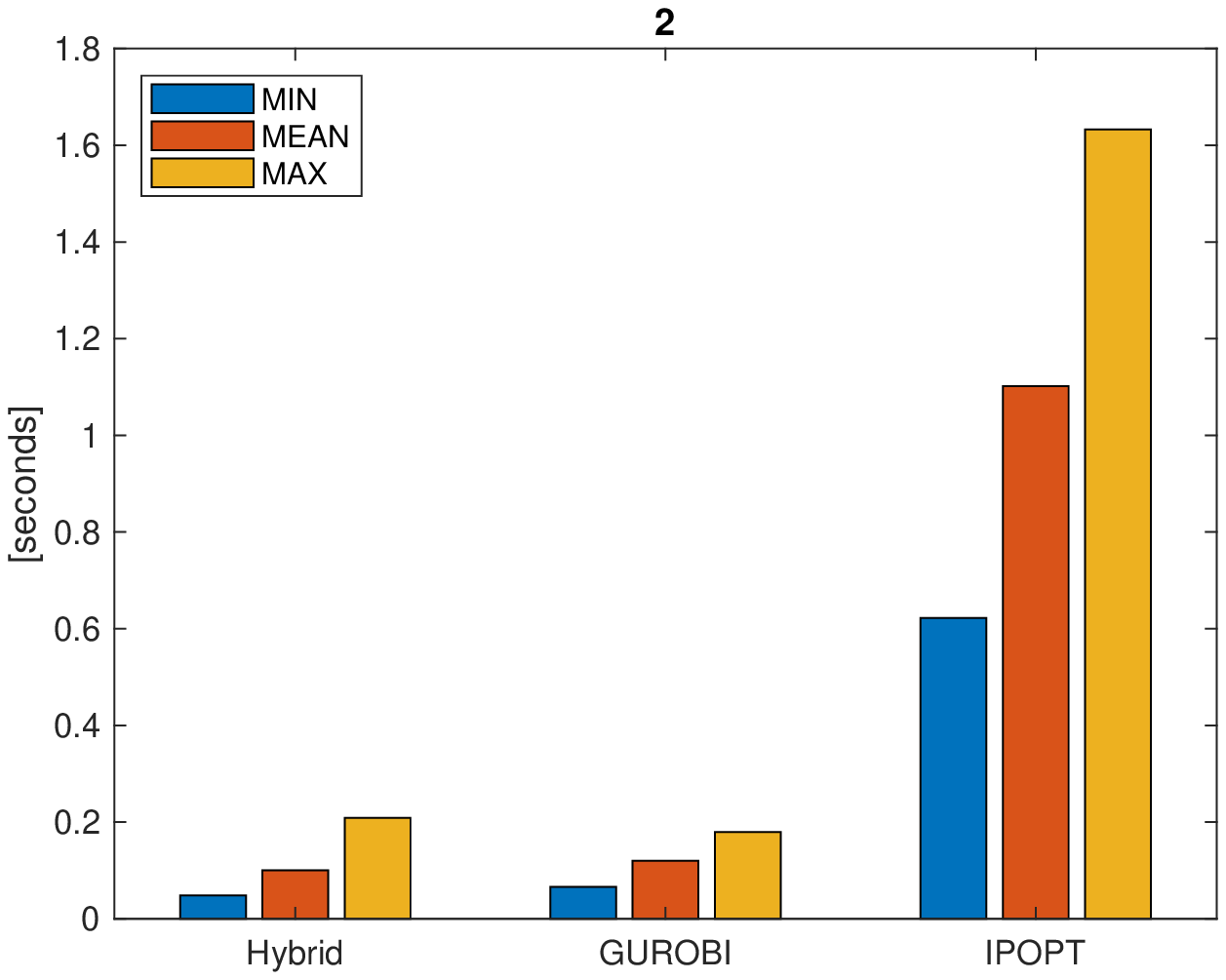}
		\caption{Samples $n = 500$}
		\label{fig:ex3_500}
	\end{subfigure}
	~ 
	\begin{subfigure}[b]{\columnwidth}
		\includegraphics[width=0.6\columnwidth]{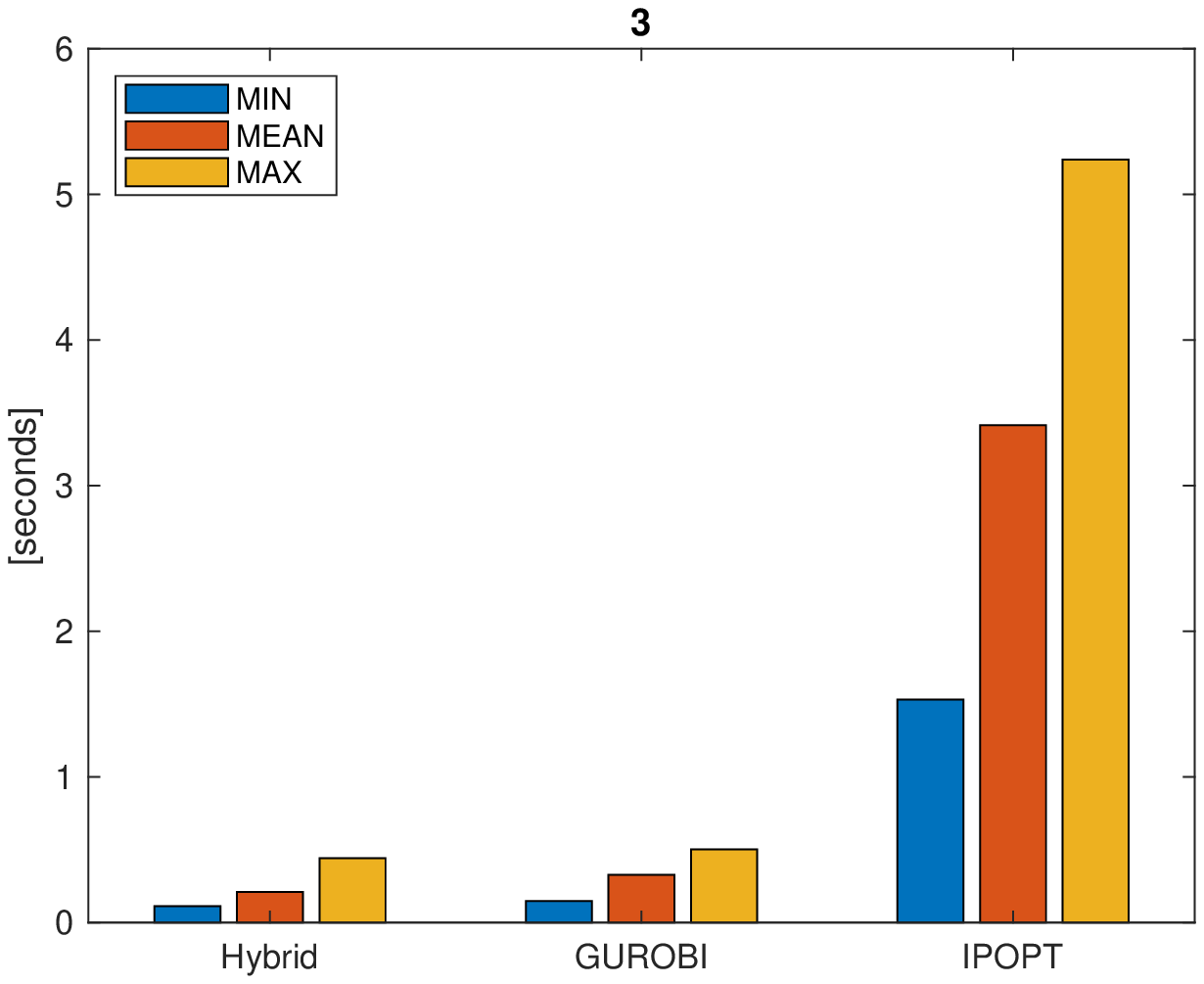}
		\caption{Samples $n = 1000$}
		\label{fig:mouse}
	\end{subfigure}
	\caption{Computational results for Experiment 2.\label{cap4:fig:ex3_result}}
\end{figure}
With respect to the paths of Experiment 1, for those in Experiment 2 the superiority of both versions of Algorithm SCA with respect to IPOPT is even clearer, and also in this case
the superiority becomes more and more evident as the number of sampled points $n$ increases. 
For these paths \emph{SCA-H} still performs quite well but it can not be claimed to be superior with respect to \emph{SCA-G}. 
Moreover, the solutions returned by  \emph{SCA-H} are in some cases poorer (in terms of objective function values) with respect to those returned by \emph{SCA-G}, although the difference is
still very mild and never exceeds 1\%. 
We can give two possible motivations:
\begin{itemize}
	\item The directions computed by the heuristic procedure are not necessarily good descent directions, so routine \texttt{computeUpdate} slowly converged to a solution. 
	\item  The heuristic procedure often failed and it was in any case necessary to call GUROBI. 
\end{itemize}
For what concerns IPOPT, besides being slower, we should also remark that for $n=100$, it is sometimes unable to converge and returns poor solutions
whose objective function values exceed by more than 100\% those returned by  \emph{SCA-H} and  \emph{SCA-G}.

\begin{figure}[!h]
	\centering
	\includegraphics[width=0.8\columnwidth]{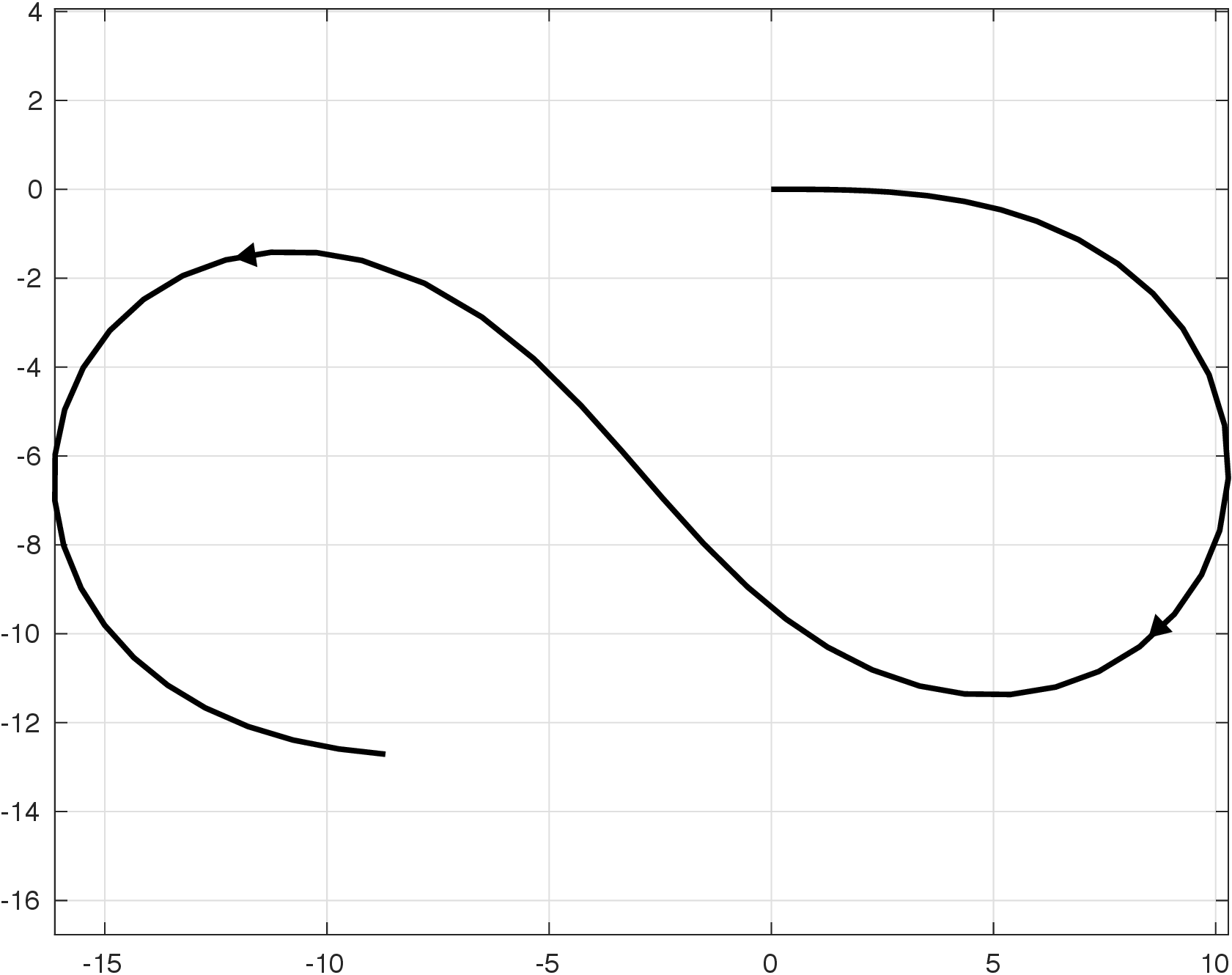}
	\caption{Experiment 3: second path.}
	\label{fig:sinpath}
\end{figure}

\begin{figure}[!h]
	\centering
	\includegraphics[width=0.8\columnwidth]{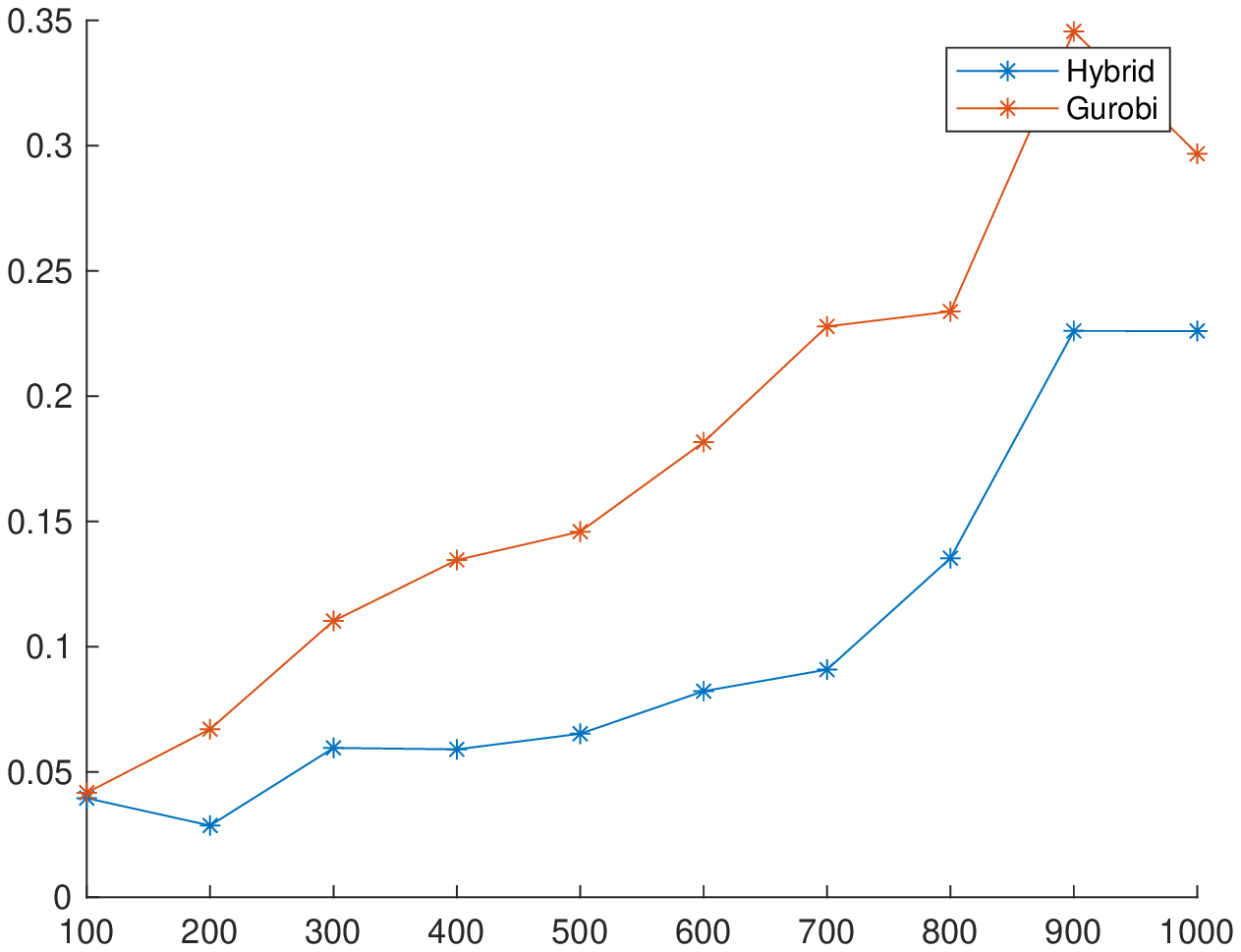}
	\caption{ Computational times as a function of the number $n$ of sample points for the two tested versions of Algorithm SCA
over the path displayed in Figure~\ref{cap1:fig:geometry_path}}
	\label{fig:minaripathtimes}
\end{figure}

\begin{figure}[!h]
	\centering
	\includegraphics[width=0.8\columnwidth]{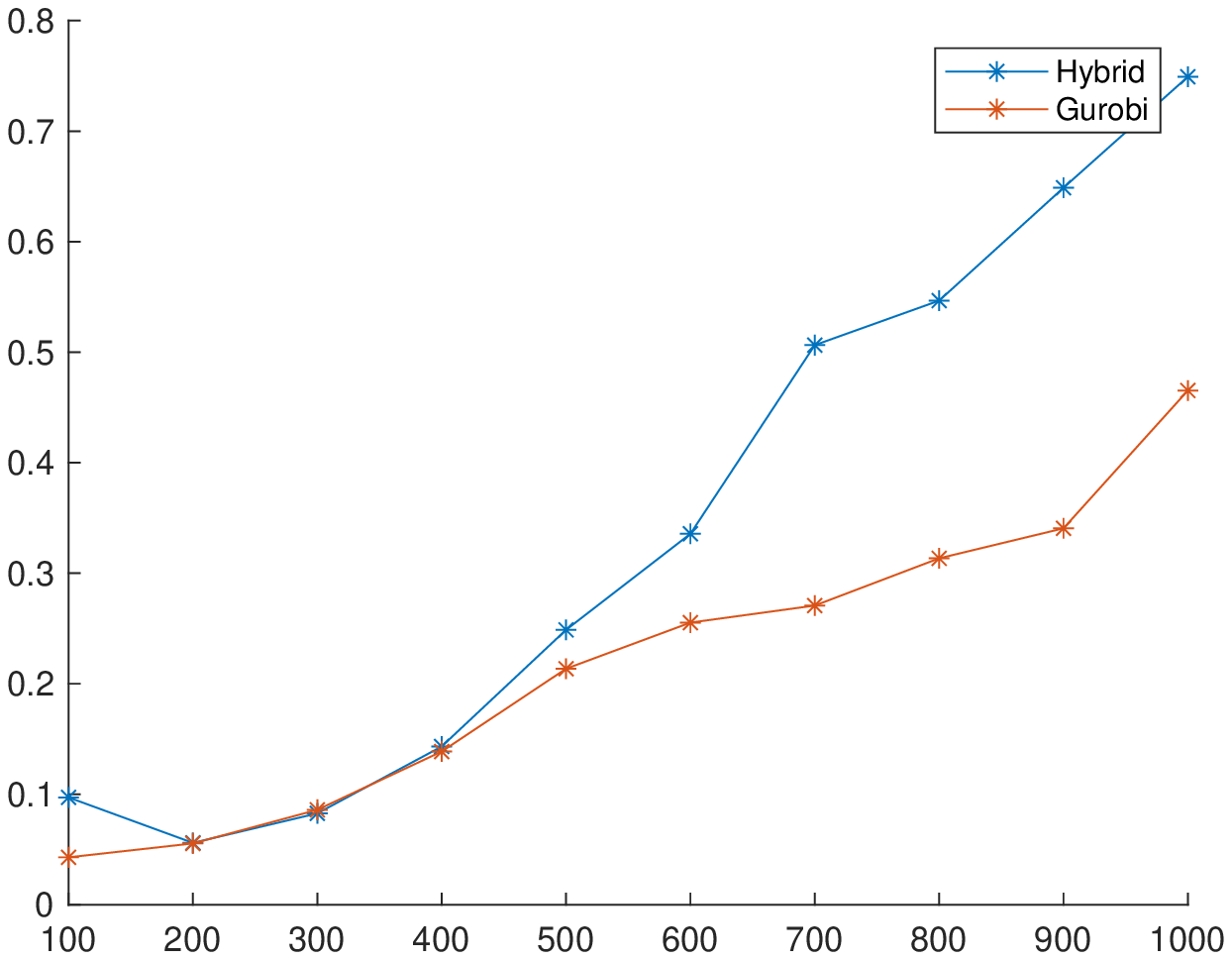}
	\caption{Computational times as a function of the number $n$ of sample points for the two tested versions of Algorithm SCA
over the path displayed in Figure~\ref{fig:sinpath}}
	\label{fig:jotatime}
\end{figure}

\paragraph{Experiment 3}
In our third experiment we compared the performance of the two proposed approaches (\emph{SCA-H} and \emph{SCA-G}) over two possible driving scenarios as the number $n$ of samples increases.
As a first example we considered a continuous curvature path composed of a line segment, a clothoid, a circle arc, a clothoid and a final line segment  (see Figure~\ref{cap1:fig:geometry_path}).
The minimum-time velocity planning on this path, whose total length is $s_{f} = 90$ m, is addressed with the following data.
The maximum squared velocity is 225 m$^2$s$^{-2}$, the longitudinal acceleration limit is $A = 1.5$ ms$^{-2}$, the maximal normal acceleration is $A_N = 1$ ms$^{-2}$, while for the jerk constraints we set $J=1$ ms$^{-3}$.
Next, we considered a path of length $s_f = 60$ m (see Figure~\ref{fig:sinpath}) whose curvature
was defined according to the following function
\[
k(s) = \frac{1}{5}\sin\left(\frac{s}{10}\right),\quad s\in[0,s_f],
\]
and parameter $A$, $A_N$ and $J$ were set equal to 1.39 ms$^{-2}$, 4.9 ms$^{-2}$ and 0.5 ms$^{-3}$, respectively.
The computational results are reported in Figure~\ref{fig:minaripathtimes} and Figure~\ref{fig:jotatime}
for values of $n$ that grows from 100 to 1000.
Figures~\ref{fig:minaripathtimes} and~\ref{fig:jotatime} show two opposite results which confirm 
what we have already observed about Experiments 1 and 2, namely, that the performance of \emph{SCA-H} and \emph{SCA-G} depend on the path. In particular, it seems that the heuristic performs in a poorer way when the number of 
points of the upper bound vector at which 
PAR constraints are violated
(which will be called critical points in Section~\ref{sec:heuristic}), tends to be large, which is the case for the second instance.
Note that, although not reported here, the computing times of IPOPT on these two paths are larger than those of \emph{SCA-H} and \emph{SCA-G}, and, as usual, the gap increases with $n$. Moreover, for the second path IPOPT
was unable to converge for $n=100$ and returned a solution which differed by more than 35\% with respect to those returned by \emph{SCA-H} and \emph{SCA-G}.

\begin{figure}[!h]
	\centering
	\includegraphics[width=1\columnwidth]{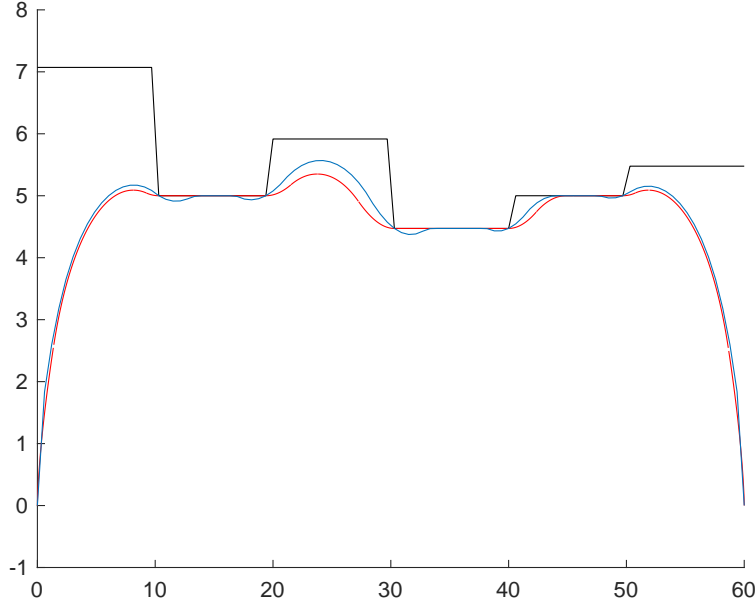}
	\caption{Velocity profile returned by the heuristic proposed in \cite{RaiCGL2019jerk} (red line) and by Algorithm SCA (blue line). The black line is the maximum velocity profile.}
	\label{fig:compareheur}
\end{figure}

\paragraph{Experiment 4} In this experiment we compared the performance of our approach with the heuristic procedure 
recently proposed in \cite{RaiCGL2019jerk}. We made different tests with the instances discussed in Experiments 1 and 2. Algorithms \emph{SCA-H} and \emph{SCA-G}  have computing times comparable (actually, slightly better)
with respect to that heuristic, and the quality of the final solutions is 5\%-10\% higher. 
Note that for a company a 5-10\% gain means the opportunity of completing 5-10\% more tasks during the day (taking into account that these algorithms are not only run once a day but are  repeatedly run throughout the day, e.g., to plan the activities of LGVs in a depot), which is a considerable gain from an economic point of view.
Rather than reporting detailed computational results, we believe that it is more instructive to discuss a single representative instance, taken from Experiment 1 with $n=100$, which reveals the qualitative difference between the solutions returned by Algorithm SCA and those returned by the heuristic. In this instance we set 
$A = 2.78$ ms$^{-2}$, while for the jerk constraints we set $J=2$ ms$^{-3}$. The total length of the path
is $s_f = 60$ m. The maximum velocity profile is the piecewise constant black line in Figure \ref{fig:compareheur}. In the same figure we report in red the velocity profile returned by the heuristic and in blue the one returned by Algorithm SCA. The computing time for the heuristic is 45ms, while for Algorithm SCA is 39ms. The final objective function value (i.e., the travelling time along the given path) is 15.4s for the velocity profile returned by the heuristic, and 14.02s for the velocity profile returned by Algorithm SCA. From the qualitative point of view it can be observed in this instance (and similar observations hold for the other instances we tested) that the heuristic produces velocity profiles whose local minima coincide with those of the maximum velocity
profile. For instance, in the interval between 10m and 20m we notice that the velocity profile returned by the heuristic coincides with the maximum velocity profile in that interval. Instead, the velocity profile generated by Algorithm SCA generates velocity profiles which fall below the
local minima of the maximum velocity profile, but this way they are able to keep the velocity higher in the regions preceding and following the local minima of the maximum velocity profile. Again  referring to the interval between 10m and 20m, we notice that the velocity profile computed by Algorithm SCA falls below the maximum velocity profile in that region and, thus, below the velocities returned by the heuristic, but this way velocities in the region before 10m and in the one after 20m are larger with respect to those computed by the heuristic.

\paragraph{Experiment 5}
As a final experiment, to illustrate the approach presented in Section~\ref{sec_higher_dim},
we consider a speed planning problem for a UAV vehicle. The
configuration space is $SE(3)$, the set of rigid
transformations in $\Real^3$. A configuration of $SE(3)$ is
represented by a couple $(R,p)$ ,with $R \in SO(3)$ (the set of rotations in
$\Real^3$) and $p \in \Real^3$. The rigid transformation associated to couple $(R,p)$ is given
by map $T: \Real^3 \to \Real^3$, 
$T(x)= R x+p$. Note that $R$ and $p$ are associated, respectively, to the vehicle rotation
and translation.
Let $(R',p')$ be an element of the tangent space of $SE(3)$ at $(R,p)$.
Then, $R'$ can be written as $R'=R \Omega$, where
$\Omega=\left[\begin{array}{lll} 0 & -\omega_3 & \omega_2 \\
          \omega_3  & 0 & -\omega_1\\
          -\omega_2 & \omega_1  & 0 \end{array}\right]$
        is a skew symmetric matrix and $\omega
        =[\omega_1,\omega_2,\omega_3]^T \in \Real^3$ is a vector that
        contains the angular velocities with respect to the vehicle
        frame.
      In our computations, we considered norm
    $\|(R',p')\|^2=\|\omega\|^2+ \|p'\|^2$. Namely, the squared
        norm of $(R',p')$ is the sum of the squared
                       angular and translational
                       velocities. Obviously, other choices are
                       possible.
                       
\begin{figure}
	\centering
     \includegraphics[width=\columnwidth]{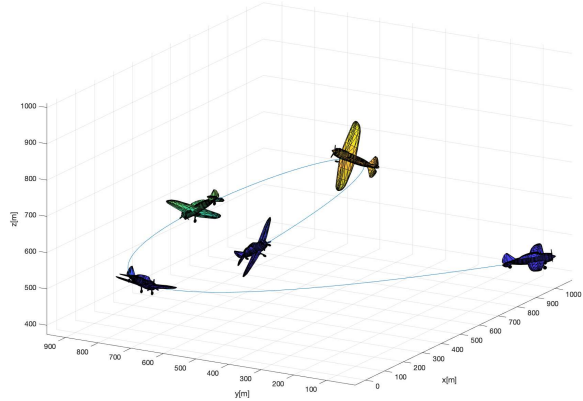}
	\caption{Example UAV reference path}
		\label{fig:path_UAV_example}
\end{figure}

\begin{figure}
	\centering
	\begin{subfigure}[b]{\columnwidth}
		\includegraphics[width=0.6\columnwidth]{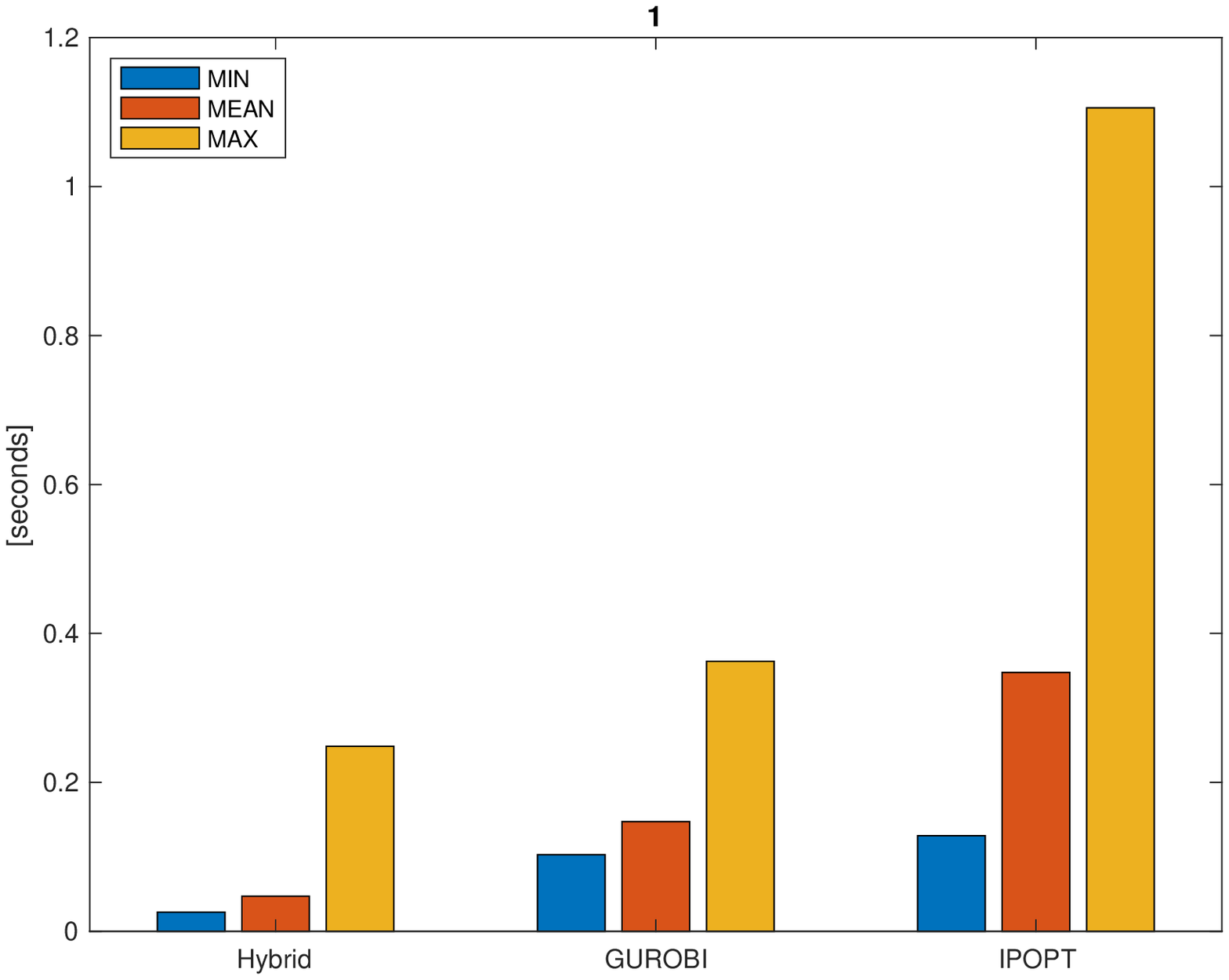}
		\caption{Samples $n = 100$}
		\label{fig:ex5_100}
	\end{subfigure}
	~ 
	\begin{subfigure}[b]{\columnwidth}
		\includegraphics[width=0.6\columnwidth]{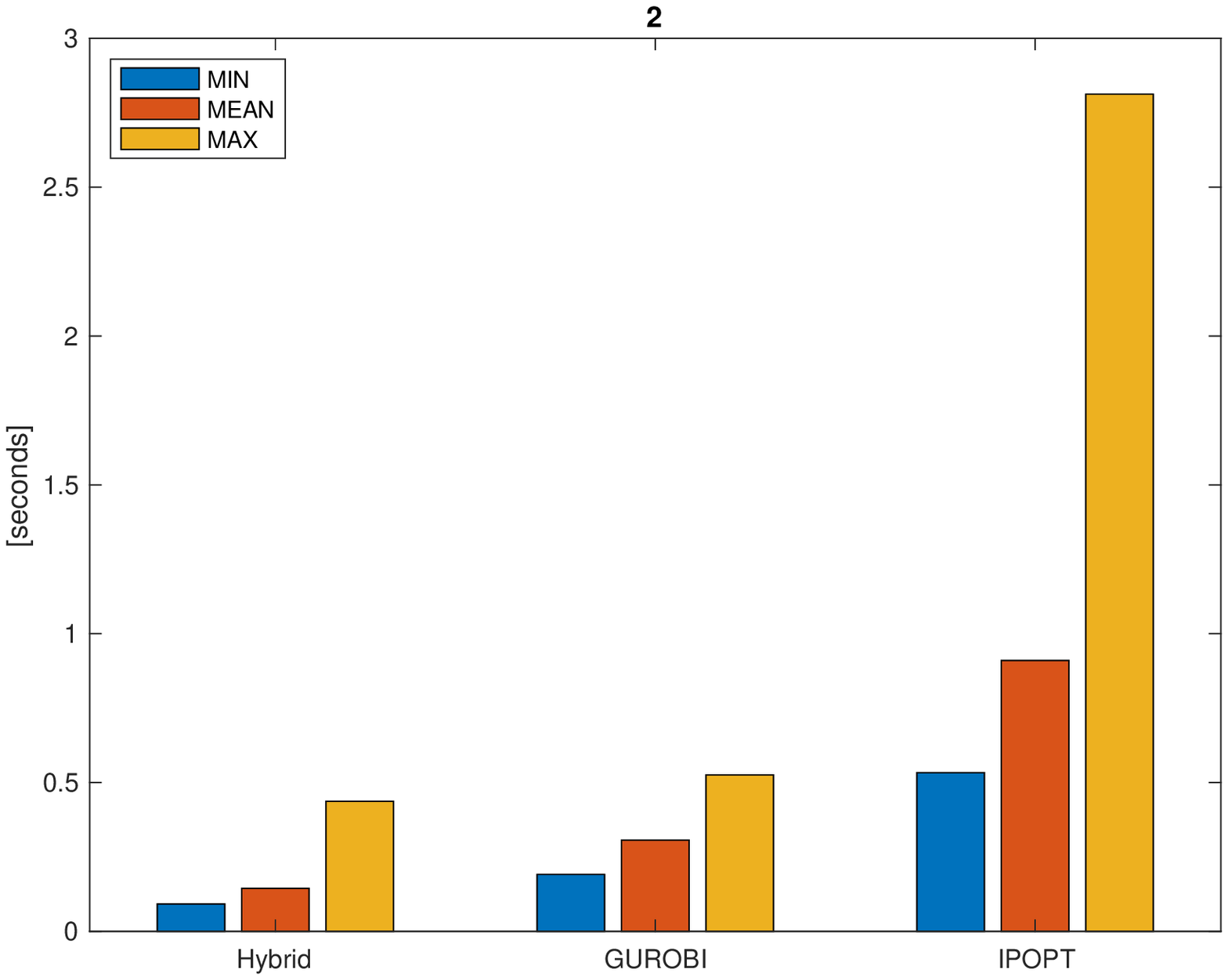}
		\caption{Samples $n = 500$}
		\label{fig: ex5_500}
	\end{subfigure}
	~ 
	\begin{subfigure}[b]{\columnwidth}
		\includegraphics[width=0.6\columnwidth]{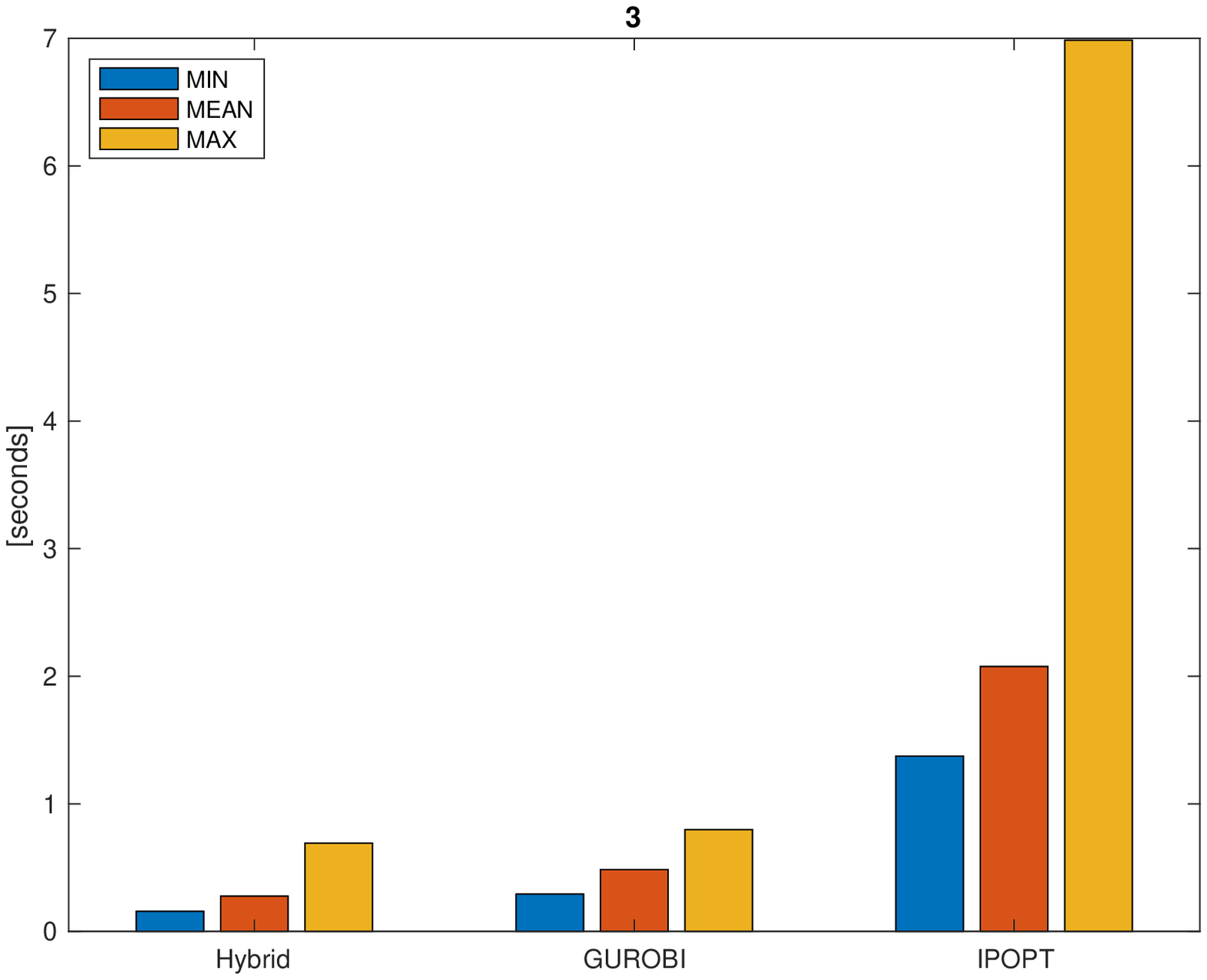}
		\caption{Samples $n = 1000$}
		\label{fig: ex5_1000}
	\end{subfigure}
	\caption{Computational results for Experiment 5.\label{cap:fig:ex5_result}}
\end{figure}

  We randomly defined a path in $SE(3)$ with the following
  procedure. We first picked independent random vectors $p_i
  \in\Real^3$, $i=0,\ldots,3$, in which each component of $p_i$ is
  chosen
  from a uniform distribution in interval $[0,1000]$. Then, we interpolated these points by defining a spline curve $\eta: [0,3]
  \to \Real^3$ of order $5$ such that $\eta(i)=p_i$, $i=0,\ldots,3$.
We associate each vector $r \in \Real^3$ to a
rotation matrix by the exponential map. That is, we define function $S:\Real^3 \to \Real^{3
  \times 3}$ such that, if $r=[x,y,z]^T$, $S(x)=\left[\begin{array}{lll}
                                                  0 &-z &y
                                                  \\z &0&-x\\-y&x&0
                                                  \end{array}\right]$
                                                    and then
                                                    set $M:
                                                    \Real^3 \to SO(3)$,
                                                    $M(x)=e^{S(x)}$.
Set $e_1=[1,0,0]^T$, the unit vector aligned with the $x$-axis.
                                                    We defined  four vectors $r_i
\in\Real^3$, $i=1\,\dots,4$, by solving equation $T(r_i) e_1=\frac{\eta'(i)}{\|\eta'(i)\|}$.
In this way, for $i=1,\ldots,4$, the $x$-axis of the vehicle frame is
aligned to the tangent of $\eta$ at $i$.
Finally,  we defined a second spline curve $\mu: [0,3]
  \to \Real^3$ of order $5$ that satisfies conditions $\mu(i)=r_i$, $i=0,\ldots,3$.
Then, the reference path is given by $\gamma: [0,3] \to SO(3) \times
R^3$, $\gamma(s)=(T(\mu(i)),\eta(i))$, after arc-length
reparameterization.
Figure~\ref{fig:path_UAV_example} presents a possible reference path
obtained with this method.
Note that this procedure is just a simple trick for determining a random
path in $SE(3)$ in order to test the procedure presented in Section~\ref{sec_higher_dim}. In general, the
determination of a reference path in $SE(3)$ for a UAV is a complex task that has to take into account multiple factors, such that
the actual dynamic model of the vehicle and the actuator
limits. However, addressing this problem is outside the scope of this
work. Indeed, the random path $\gamma$ obtained with the method
presented here may not be a valid reference for a UAV. 
Figures~\ref{fig:ex5_100}--\ref{fig: ex5_1000} show the computation
times for algorithms \emph{SCA-H}, \emph{SCA-G} and IPOPT, for $n \in
\{100,500,1000\}$. We applied the method presented in
Section~\ref{sec_higher_dim} with $\hat V=50$, $\hat A=5$, $\hat J=1$,
$A=\frac{\hat A}{2}$, $\hat J=\frac{\hat J}{2}$.
The results are qualitatively similar to previous experiments.
                                           
\bibliographystyle{unsrt}    
\bibliography{VelPlan,Automatica,Automatica_PR} 
\begin{appendix}
\section{Proof of Theorem \ref{thm:convergence}}
\label{sec:appconvergence}
In order to prove the theorem, we first need to prove some lemmas.
\begin{lem}
\label{lem:convfun}
The sequence
$\{f({\bf w}^{(k)})\}$ of the function values at points generated by Algorithm SCA converges to a finite value.
\end{lem}
\begin{proof}
The sequence is nonincreasing and bounded from below, e.g., by the value $f({\bf u}_B)$, in view of the fact that the objective function $f$ is monotonic decreasing. Thus, it converges to a finite value.
\end{proof}
Next, we need the following result based on strict convexity of the objective function $f$.
\begin{lem}
\label{lem:strconv}
For each $\delta>0$ sufficiently small, it holds that
\begin{equation}
\label{eq:auxstrconv}
\begin{array}{ll}
\min & \left\{\max\{f({\bf x}), f({\bf y})\}-f\left(\frac{{\bf x}+{\bf y}}{2}\right)\ :\ \right. \\ [6pt]
       & \left.{\bf x}, {\bf y}\in \Omega, \ \|{\bf x}-{\bf y}\|\geq \delta\right\}\geq \varepsilon_{\delta}>0.
\end{array}
\end{equation}
\end{lem}
\begin{proof}
Due to strict convexity, it holds that $\forall {\bf x}\neq {\bf y}$,
$$
\max\{f({\bf x}), f({\bf y})\}-f\left(\frac{{\bf x}+{\bf y}}{2}\right)>0.
$$
Moreover, the function is a continuous one.
Next, we observe that the region 
$$
\left\{{\bf x}, {\bf y}\in \Omega:\ \|{\bf x}-{\bf y}\|\geq \delta\right\},
$$
is a compact set. Thus, by Weierstrass Theorem, the minimum in (\ref{eq:auxstrconv}) is attained and it must be strictly positive, as we wanted to prove.
\end{proof}
Finally, we prove that also the sequence of points generated by Algorithm SCA converges to some point, feasible for 
Problem \ref{cap4:prob_disc}.
\begin{lem}
\label{lem:seqpoint}
It holds that
$$
\|\delta {\bf w}^{(k)}\|\rightarrow 0.
$$
\end{lem}
\begin{proof}
Let us assume, by contradiction, that over some infinite subsequence with index set ${\cal K}$,
it holds that $\|\delta {\bf w}^{(k)}\|\geq 2\rho>0$ for all $k\in {\cal K}$, i.e.,
\begin{equation}
\label{eq:newpoint}
\|{\bf w}^{(k+1)}-{\bf w}^{(k)}\|\geq 2\rho>0,
\end{equation}
where ${\bf w}^{(k+1)}={\bf w}^{(k)}+\delta {\bf w}^{(k)}$.
Over this subsequence it holds, by strict convexity, that 
\begin{equation}
\label{eq:strdecr}
f({\bf w}^{(k+1)})\leq f({\bf w}^{(k)})-\xi\  \ \ \ \ \forall k\in {\cal K},
\end{equation} 
for some $\xi>0$. Indeed, it follows by optimality of ${\bf w}^{(k)}+\delta {\bf w}^{(k)}$ for Problem \ref{cap4:prob_lin} and convexity of $f$ that
$$
f({\bf w}^{(k+1)})\leq f\left(\frac{{\bf w}^{(k+1)}+{\bf w}^{(k)}}{2}\right)\leq f({\bf w}^{(k)}),
$$
so that
$$
\max\left\{f({\bf w}^{(k)}), f({\bf w}^{(k+1)}\right\}=f({\bf w}^{(k)}).
$$
Then, it follows from  (\ref{eq:newpoint}) and Lemma 
\ref{lem:strconv} that we can choose $\xi=\varepsilon_{\rho}>0$.
Thus, since (\ref{eq:strdecr}) holds infinitely often, we should have
$f({\bf w}^{(k)})\rightarrow -\infty$, which, however, is not possible in view of Lemma \ref{lem:convfun}.
\end{proof}
As a consequence of Lemma \ref{lem:seqpoint} it also holds that
\begin{equation}
\label{eq:convergpoint}
{\bf w}^{(k)}\rightarrow \bar{{\bf w}}\in \Omega.
\end{equation}
Indeed, all points ${\bf w}^{(k)}$ belong to the compact feasible region $\Omega$, so that the sequence
$\{{\bf w}^{(k)}\}$ admits accumulation points. However, due to Lemma \ref{lem:seqpoint}, the sequence cannot have
distinct accumulation points.
\newline\newline\noindent
Now, let us consider the compact reformulation (\ref{eq:origcompact}) of Problem  \ref{cap4:prob_disc} and the related
linearization (\ref{eq:lincompact}), equivalent to Problem  \ref{cap4:prob_lin} with the linearized constraints
\eqref{cap4:con:linnar1}-\eqref{cap4:con:linpar1}. 
Since the latter is a convex problem with linear constraints, its local minimizer $\delta {\bf w}^{(k)}$ (unique in view of strict convexity of the objective function) fulfills the following KKT conditions
\begin{equation}
\label{eq:sysKKTlin}
\begin{array}{l}
\nabla f({\bf w}^{(k)}+\delta {\bf w}^{(k)}) +\boldsymbol{\mu}_k^\top \nabla {\bf c}({\bf w}^{(k)})={\bf 0} \\ [6pt]
{\bf c}({\bf w}^{(k)}) +\nabla {\bf c}({\bf w}^{(k)}) \delta {\bf w}^{(k)} \leq 0 \\ [6pt]
\boldsymbol{\mu}_k^\top\left({\bf c}({\bf w}^{(k)}) +\nabla {\bf c}({\bf w}^{(k)}) \delta {\bf w}^{(k)} \right)=0 \\ [6pt]
\boldsymbol{\mu}_k\geq {\bf 0},
\end{array}
\end{equation}
where $\boldsymbol{\mu}_k$ is the vector of Lagrange multipliers.
Now, by taking the limit of system (\ref{eq:sysKKTlin}), possibly over a subsequence, in order to guarantee convergence of the multiplier vectors $\boldsymbol{\mu}_k$
to a vector $\bar{\boldsymbol{\mu}}$, in view of Lemma \ref{lem:seqpoint}
and of (\ref{eq:convergpoint}), we have that
$$
\begin{array}{l}
\nabla f(\bar{{\bf w}}) +\bar{\boldsymbol{\mu}}^\top \nabla {\bf c}(\bar{{\bf w}})={\bf 0} \\ [6pt]
{\bf c}(\bar{{\bf w}})  \leq 0 \\ [6pt]
\bar{\boldsymbol{\mu}}^\top {\bf c}(\bar{{\bf w}})=0 \\ [6pt]
\bar{\boldsymbol{\mu}}\geq{\bf  0},
\end{array}
$$
or, equivalently, the limit point $\bar{{\bf w}}$ is a KKT point of Problem \ref{cap4:prob_disc}, as we wanted to prove.

\section{Proof of Proposition \ref{cap4:prop:salvezza}}
\label{sec:feasdir}
First, we notice that 
if we prove the result for the tighter constraints \eqref{cap4:con:linnar1}-\eqref{cap4:con:linpar1}, then it must also hold
for constraints \eqref{cap4:eq:nar_lin}-\eqref{cap4:eq:par_lin}. So we prove the result only for the former.
By definition (\ref{eq:computedir}), $\dwb$ satisfies the acceleration and NAR constraints, so that 	
	$$
	\begin{array}{l}
	\dw_j\leq \dw_{j+1}+ b_{D_j} \\ [8pt]
	\dw_j\leq \dw_{j-1}+ b_{A_j} \\ [8pt]
	\dw_j \leq \beta_j(\dw_{j+1}+\dw_{j-1})+ b_{N_j} \\ [8pt]
	\dw_j\leq y^*_j.
	\end{array}
	$$
	At least one of these constraints must be active, otherwise $\dw_j$ could be increased, thus contradicting optimality.
	If the active constraint is $	 \dw_j \leq \beta_j(\dw_{j+1}+\dw_{j-1})+ b_{N_j}$, then constraint \eqref{cap4:con:linpar1} can be rewritten as follows
$$
(\theta_j-\beta_j)(\dw_{j+1}+\dw_{j-1})\leq b_{P_j}+b_{N_j}.
$$	
By recalling the definitions of $\theta_j, \beta_j,  b_{P_j}$, and $b_{N_j}$, it can be seen that this is equivalent to (\ref{cap4:ass:prop_salvezza}) and, thus, the constraint is satisfied under the given assumption.
	If $\dw_j = y^*_j$, then 
	$$
	\theta_j(\dw_{j-1}+\dw_{j+1})\leq\theta_j( y^*_{j-1}+y^*_{j+1})\leq  y^*_j +  b_{P_j} =  \dw_j +  b_{P_j} ,
	$$
	where the second inequality follows from the fact that ${\bf y}^*$ satisfies the PAR constraints. 
	Now, let $\dw_j =  \dw_{j+1}+ b_{D_j}$ (the case when $\dw_j\leq \dw_{j-1}+ b_{A_j}$ is active can be dealt with in a completely analogous way).
	First we observe that   $\dw_j\geq \dw_{j-1}- b_{D_{j-1}}$. 
	Then, 
	\begin{equation*}
	2\dw_j \ge \dw_{j+1} + \dw_{j-1} + b_{D_j} - b_{D_{j-1}}.
	\end{equation*}
	In view of the definitions of $b_{D_j}$  and $ b_{D_{j-1}}$ this can also be written as 
		\begin{equation}
\label{cap4:aux_ineq}
	2\dw_j \ge \dw_{j+1} + \dw_{j-1} + w^{(k)}_{j+1} -2 w^{(k)}_{j}+ w^{(k)}_{j-1}.
	\end{equation}
	Now we recall that 
$$
\begin{array}{l}
	\theta_j = \frac{1}{2} + \frac{\Delta}{2\left(w^{(k)}_{j+1} + w^{(k)}_{j-1}\right)^{\frac{3}{2}}} \\ [6pt] 
b_{P_j}=\frac{\Delta}{\left(w^{(k)}_{j+1} + w^{(k)}_{j-1}\right)^{\frac{1}{2}}} - \frac{1}{2}\left(w^{(k)}_{j+1} -2 w^{(k)}_{j}+ w^{(k)}_{j-1}\right),
\end{array}
$$
where $\Delta= \sqrt{2}h^2J$. Then,
\eqref{cap4:con:linpar1} can be rewritten as
	\[
	2\dw_j \ge \dw_{j+1} + \dw_{j-1} + \frac{\Delta}{\left(w^{(k)}_{j+1} + w^{(k)}_{j-1}\right)^{\frac{3}{2}}} (\dw_{j+1} + \dw_{j-1})-2b_{P_j}.
	\] 
	Now, taking into account~\eqref{cap4:aux_ineq}, such inequality certainly
	holds if 
$$
w^{(k)}_{j+1} -2 w^{(k)}_{j}+ w^{(k)}_{j-1}\geq \frac{\Delta}{\left(w^{(k)}_{j+1} + w^{(k)}_{j-1}\right)^{\frac{3}{2}}} (\dw_{j+1} + \dw_{j-1})-2b_{P_j}.
$$
Recalling the definition of $b_{P_j}$, the above inequality can be rewritten as
	\[
	\frac{2\Delta}{\sqrt{w^{(k)}_{j+1} + w^{(k)}_{j-1}}}\ge \frac{\Delta}{\left(w^{(k)}_{j+1} + w^{(k)}_{j-1}\right)^{\frac{3}{2}}} (\dw_{j+1} + \dw_{j-1}),
	\]
	and it holds if $2\left(w^{(k)}_{j+1} + w^{(k)}_{j-1}\right) \ge(\dw_{j+1} + \dw_{j-1}) $, as we wanted to prove.
\section{A heuristic procedure for computing a suboptimal descent direction\label{cap4:subsec:heur}}
\label{sec:heuristic}
We first need to introduce a definition.
\begin{defn}
	Given a vector $\bd\in\Real^N$ the set of critical points associated to such vector is
	\[
	Q(\bd) = \left\{ p \ : \  \eta_i (d_{p-1}  + d_{p+1}) -d_{p} - \beta_p > 0  \right\},
	\]
	i.e., the set of points where constraints~(\ref{cap4:con:par_dir_tr}) are violated at $\bd$.
\end{defn}
Now, the heuristic is detailed in Algorithm~\ref{cap4:alg:heur}.  
Its purpose is to sequentially remove all the critical points $p$ of the upper bound ${\mathbf{\bar{u}_B}}$ by activating a sequence of constraints~\eqref{cap4:con:par_dir_tr} in the neighbourhood of $p$ itself.  
We initially set $\bd = \mathbf{\bar{u}_{B}}$ and
compute the related set $Q(\mathbf{\bar{u}_{B}})$ of critical points. After that, we consider the most
violated critical point $p\in Q(\mathbf{\bar{u}_{B}})$ and define $\Delta_p$ as its  associated violation.
Then, we define the \emph{propagation function} from $p$. To this aim we first define a function
${\bf z}: [0,1] \rightarrow \Real^n$ such that:
\begin{equation}\label{cap4:eq:propag}
z_j(\alpha; \bd,p) = \begin{cases}
d_p & j=p\\
d_{p-1} - \alpha\Delta_p & j=p-1\\
d_{p+1} - (1-\alpha)\Delta_p & j=p+1\\
\eta_{j}^{-1}(\beta_{j} + z_{j}(\alpha; \bd,p)) - z_{j+1}(\alpha; \bd,p), & j <p \\ 
\eta_{j}^{-1}(\beta_{j} + z_{j}(\alpha; \bd,p)) - z_{j-1}(\alpha; \bd,p), & j >p. \\ 
\end{cases}
\end{equation}
Then, let 
$$
\begin{array}{l}
k_1=\max\{k<p\ :\ z_k(\alpha; \bd,p)\geq d_k\} \\ [8pt]
k_2=\min\{k>p\ :\ z_k(\alpha; \bd,p)\geq d_k\}.
\end{array}
$$
We define the propagation function ${\bf x}(\cdot; \bd, p): [0,1] \rightarrow \Real^n$ around $p$ as follows:
\begin{equation}\label{cap4:eq:propag}
x_j(\alpha; \bd; p) = \begin{cases}
z_j(\alpha; \bd,p) & j=k_1+1,\ldots,k_2-1\\
d_j & \mbox{otherwise.}
\end{cases}
\end{equation}
Basically, ${\bf x}$ decreases the components of the current vector $\bd$ around the critical point $p$ in order to remove the violations locally, without decreasing all other components.
The choice of not decreasing the remaining components comes from the fact that the objective function \eqref{cap4:obj:d} is monotonic non-increasing.
After having defined the propagation function, we search for the best $\alpha$ that, in addition to activating a sequence of constraints~\eqref{cap4:con:par_dir_tr} around $p$, gives the best possible solution with respect to the objective function~\eqref{cap4:obj:d}.
Then, we consider the following problem:
\begin{equation}\label{cap4:prob:under-estimator}
\alpha^*\in \arg \min_{\alpha\in[0,1]} \left\{ -\nub^T{\xb}(\alpha), | , {\xb}(\alpha) \ge \mathbf{\bar{l}_B} \right\}.
\end{equation}
We employ the ternary search algorithm
to efficiently find its optimal solution~\cite{cormen2009introduction}. 
Note that the ternary search needs to consider the lower bound constraints.
Actually, this issue can be easily overcome by setting to $+\infty$ the value of the objective function if  ${\bf x}(\alpha) \not\geq {\bf {\bar{l}}_{B}}$.
If an optimal solution $\alpha^*$ exists, we set $\bd = \xb(\alpha^*)$, compute its set of critical points $Q(\bd)$, and repeat the above procedure until $Q$ is empty.
If an optimal solution $\alpha^*$ does not exist, i.e., problem~\eqref{cap4:prob:under-estimator} is unfeasible, then
we remove the critical point $p$ from $Q$ and repeat the above procedure by considering the next most violated critical point.
\begin{remark}\label{cap4:rem:heur}
	The procedure described in this section is a heuristic one since we do not have any proof of correctness and optimality.	
	Moreover, it may happen that:
	\begin{itemize}
		\item $\alpha^*$ does not exist for all critical points contained in the set $Q$, i.e., we are unable to remove all violations;
		\item the solution returned by the heuristic might not be a descent direction, i.e., $-\nub^T\bd \ge 0$.
	\end{itemize}
	For this reason, if one of these two cases occurs, we do not consider the result computed by the heuristic and solve Problem~\ref{cap4:prob:direction_tr} by using an LP solver.
\end{remark}
\begin{algorithm}[!h]
	\caption{The heuristic procedure to compute a descent direction ${\bd}$ for Problem~\ref{cap4:prob:direction_tr}\label{cap4:alg:heur}}
	\SetKwProg{function}{Function}{:}{}
	Set $\bd = \mathbf{\bar{u}_B}$\;
	Let $U=Q({\bd})$\;
	\While{$U \ne \emptyset$}{
		Let $p$ be the most violated critical point in $U$\;
		Let ${\bf x}(\alpha; \bd; p)$ be defined as in (\ref{cap4:eq:propag})\;
		\If{$ \exists \ \alpha^* \in \arg\min_{\alpha\in[0,1]} \left\{ -\nub^T \xb(\alpha; \bd, p) \, | \, \xb(\alpha; \bd,p) \ge \mathbf{\bar{l}_B} \right\}$}
		{	Set $\bd= \xb(\alpha^*; \bd, p)$\;
			Let $U=Q(\bd)$
		}
		\Else{
			$U = U \backslash \{p\}$ 
		}
	}
	
	\Return $\bd$
	%
	%
\end{algorithm}


\end{appendix}
\end{document}